\documentclass[10pt,twoside,reqno]{amsart}

\usepackage{amssymb}
\usepackage{amsmath}
\usepackage{mathrsfs}
\usepackage{dsfont}
\usepackage{epsfig}
\usepackage{epstopdf}
\usepackage{url}
\usepackage{booktabs}
\usepackage{caption}

\usepackage{float}
\floatstyle{plain}
\restylefloat{figure}

\usepackage[math]{cellspace}
\newtheorem{theorem}{\sc\hspace{8pt} Theorem}
\newtheorem{corollary}[theorem]{\sc\hspace{8pt} Corollary}

\newcommand{\abs}[1]{\lvert#1\rvert}

\DeclareMathOperator*{\im}{Im}
\DeclareMathOperator*{\re}{Re}

\makeatletter
\renewcommand*\env@matrix[1][c]{\hskip -\arraycolsep
  \let\@ifnextchar\new@ifnextchar
  \array{*\c@MaxMatrixCols #1}}
\makeatother

\makeatletter
\newcommand{\vast}{\bBigg@{3.5}}
\newcommand{\Vast}{\bBigg@{4}}
\makeatother

\begin{document}

\title[The density of complex zeros of random sums]{The density of complex zeros of random sums}

\author[Christopher Corley]{Christopher Corley}

\address{Department of Mathematics, University of Tennessee at Chattanooga, 415 EMCS Building (Dept. 6956), 615 McCallie Avenue, Chattanooga, Tennessee 37403, United States of America}

\email{Christopher-Corley@utc.edu}

\author[Andrew Ledoan]{Andrew Ledoan$^{\ast}$}

\address{Department of Mathematics, University of Tennessee at Chattanooga, 415 EMCS Building (Dept. 6956), 615 McCallie Avenue, Chattanooga, Tennessee 37403, United States of America}

\email{Andrew-Ledoan@utc.edu}

\thanks{$^{\ast}$The research of the second author was supported by the National Science Foundation under Grant DMS-1852288 and by the Center of Excellence in Applied Computational Science and Engineering under Grant FY2019 CEACSE Awards.}

\subjclass[2010]{Primary 30C15; Secondary 30B20, 60B99, 60E05, 60G50, 62H10}

\keywords{Expected number of complex zeros; density of zeros; multivariate Gaussian distribution; level crossing; random polynomial; random sum; Adler's theorem; Christoffel--Darboux formula; orthogonal polynomials}

\begin{abstract}
Let $\{\eta_{j}\}_{j = 0}^{N}$ be a sequence of independent, identically distributed random complex Gaussian variables, and let $\{f_{j} (z)\}_{j = 0}^{N}$ be a sequence of given analytic functions that are real-valued on the real number line. We prove an exact formula for the expected density of the distribution of complex zeros of the random equation $\sum_{j = 0}^{N} \eta_{j} f_{j} (z) = \mathbf{K}$, where $\mathbf{K} \in \mathds{C}$. The method of proof employs a formula for the expected absolute value of quadratic forms of Gaussian random variables. We then obtain the limiting behaviour of the density function as $N$ tends to infinity and provide numerical computations for the density function and empirical distributions for random sums with certain functions $f_{j} (z)$. Finally, we study the case when the $f_{j} (z)$ are polynomials orthogonal on the real line and the unit circle.
\end{abstract}

\maketitle

\thispagestyle{empty}

\section{Introduction}

Let $\{a_{j}\}_{j = 0}^{N}$ and $\{b_{j}\}_{j = 0}^{N}$ be sequences of mutually independent, identically distributed random real variables defined on the complete probability space $(\Omega, \mathscr{F}, P)$, each sequence normally distributed with mean zero and variance one. As usual, $\Omega$ is a set with generic elements $\omega$, $\mathscr{F}$ is a $\sigma$-field of subsets of $\Omega$, and $P$ is a probability measure on $\mathscr{F}$. Throughout this paper, we shall assume that all sub $\sigma$-fields contain all sets of measure zero. Let $\{f_{j} (z)\}_{j = 0}^{N}$ be a sequence of given analytic functions that are real-valued on the real number line. Furthermore, let us define the random sum
\begin{equation*}
S_{N} (z)
 = \sum_{j = 0}^{N} \eta_{j} f_{j} (z)
 = \eta_{0} f_{0} (z) + \ldots + \eta_{N} f_{N} (z),
\end{equation*}
where $z$ is the complex variable $x + i y$, and the $\eta_{j}$ are independent, identically distributed random complex Gaussian variables (with density $e^{- z \bar{z}} / \pi$) given by $\eta_{j} = a_{j} + i b_{j}$. We denote by $\nu_{N, \mathbf{K}} (\Phi)$ the number of complex zeros in a compact subset $\Phi \subset \mathds{C}$ of $S_{N} (z)$ with respect to the complex level $\mathbf{K} = (K_{1}, K_{2})^{\prime}$. We do not assume necessarily that the scalars $K_{1}$ and $K_{2}$ are equal. From \cite{SheppVanderbei1995}, we see that, with probability one, the value of the density function $h_{N, \mathbf{K}} (z)$ for multivariate Gaussian coefficients is given by
\begin{equation} \label{equation-1}
\mathscr{E} \nu_{N, \mathbf{K}} (\Phi)
 = \int_{\Phi} h_{N, \mathbf{K}} (z) \, d z,
\end{equation}
where $\mathscr{E} \nu_{N, \mathbf{K}} (\Phi)$ is the mathematical expectation of $\nu_{N, \mathbf{K}} (\Phi)$. Thus, $h_{N, \mathbf{K}} (z)$ is the expected density of the complex zeros of the random equation $S_{N} (z) = \mathbf{K}$.

Shepp and Vanderbei \cite{SheppVanderbei1995} wrote a beautiful paper on the complex zeros of the random polynomial $\sum_{j = 0}^{N} \eta_{j} z^{j}$, where the $\eta_{j}$ are independent, identically distributed random real Gaussian coefficients. In their paper, Shepp and Vanderbei introduced a sophisticated method based on Cauchy's argument principle for producing an explicit density function $h_{N, \mathbf{0}} (z)$ for the zeros. The method uses the Cholesky decomposition for representing correlated random Gaussian variables in terms of uncorrelated (and hence independent) random Gaussian variables. Shepp and Vanderbei generated computer plots of this density function and hundreds of thousands of zeros from randomly generated polynomials that show that, as the degrees $N$ become large, the zeros tend to lie very close to the unit circle and, when the real zeros are ignored, appear to be approximately uniformly distributed around the unit circle. Their asymptotics for the density function confirm the classical result due to Hammersley \cite{Hammersley1956}.

Ibragimov and Zeitouni \cite{IbragimovZeitouni1997} obtained the results in \cite{SheppVanderbei1995} using a different method, based on an integral representation of the average number of zeros of a random field. Furthermore, Ibragimov and Zeitouni demonstrated the limiting distributions for $h_{N, \mathbf{0}} (z)$ under more general distributional assumptions.

In later work, Vanderbei \cite{Vanderbei2015} introduced a modest generalization to the central assumptions underlying the results in \cite{SheppVanderbei1995} and showed that comparable explicit formulas for the distribution of complex zeros can still be obtained for any $N$. Following the same general methodology given in \cite{SheppVanderbei1995}, Vanderbei derived analogous explicit formulas for the density of complex zeros of the random sum $S_{N} (z)$, for which the $\eta_{j}$ are independent, identically distributed random real Gaussian coefficients, and the $f_{j} (z)$ are given analytic functions that are real-valued on the real number line.

In this paper, we shall continue the line of investigation begun by Vanderbei and study the number of level crossings of the random sum $S_{N} (z)$. To consider this challenging general case, we shall employ a multivariate analysis approach based on results due to Adler \cite{Adler1981}, which provide a representation for the expected number of zeros of certain random fields. The method of proof was first applied by Ibragimov and Zeitouni \cite{IbragimovZeitouni1997}. Our main result generalizes the density functions obtained by Shepp and Vanderbei \cite{SheppVanderbei1995}, and later by Vanderbei \cite{Vanderbei2015}, to nonzero $\mathbf{K}$.

\begin{theorem} \label{theorem-1}
Let the density function $h_{N, \mathbf{K}} (z)$ be defined by \eqref{equation-1}. Under the conditions imposed on the random sum $S_{N} (z)$ and $\{\eta_{j}\}_{j = 0}^{N}$, for all integers $N > 0$ we have
\begin{equation*}
\begin{split}
h_{N, \mathbf{K}} (z)
 &= \frac{e^{-(K_{1}^{2} + K_{2}^{2}) / 2 B_{0, N} (z)}}{\pi B_{0, N} (z)} \left\{B_{2, N} (z) - \frac{\abs{B_{1, N} (z)}^{2}}{B_{0, N} (z)} \left(1 - \frac{K_{1}^{2} + K_{2}^{2}}{2 B_{0, N} (z)}\right)\right\},
\end{split}
\end{equation*}
where the kernels $B_{r, N} (z)$ for $0 \leq r \leq 2$ are given by
\begin{equation*}
\begin{array}{lll}
\displaystyle B_{0, N} (z)
 = \sum_{j = 0}^{N} \abs{f_{j} (z)}^{2},
& \displaystyle B_{1, N} (z)
 = \sum_{j = 0}^{N} \overline{f_{j} (z)} f^{\prime}_{j} (z),
& \displaystyle B_{2, N} (z)
 = \sum_{j = 0}^{N} \abs{f^{\prime}_{j}(z)}^{2}.
\end{array}
\end{equation*}
\end{theorem}

The value of $h_{N, \mathbf{K}} (z)$ is expressed in fairly simple terms, in that we can clearly see its form of dependence on $\mathbf{K}$. Its form is reminiscent of Farahmand and Jahangiri's \cite[Theorem 1]{FarahmandJahangiri1998} expected density for the random polynomial $\sum_{j = 0}^{N} \eta_{j} g_{j} z^{j}$ with respect to $\mathbf{K}$, where the $g_{j}$ are given real constants. (See \cite{Farahmand1997, Farahmand1998} for the case when $g_{j} = 1$.) From Theorem \ref{theorem-1}, we have the following consequence.

\begin{corollary} \label{corollary-1}
For any vector $\mathbf{K}$ restricted to a circle of radius $K > 0$ and for all integers $N > 0$ we have
\begin{equation*}
h_{N, \mathbf{K}} (z)
 = \frac{e^{-K^{2} / B_{0, N} (z)}}{\pi B_{0, N} (z)} \left\{B_{2, N} (z) - \frac{\abs{B_{1, N} (z)}^{2}}{B_{0, N} (z)} \left(1 - \frac{K^{2}}{B_{0, N} (z)}\right)\right\}.
\end{equation*}
\end{corollary}

Immediate by Corollary \ref{corollary-1} is the following consequence, which was proved independently by Yeager \cite{Yeager2016} and one of the authors \cite{Ledoan2016}.

\begin{corollary} \label{corollary-2}
For any vector $\mathbf{K}$ of equal elements $K = 0$ and for all integers $N > 0$ we have
\begin{equation*}
h_{N, \mathbf{0}} (z)
 = \frac{B_{0, N} (z) B_{2, N} (z) - \abs{B_{1, N} (z)}^{2}}{\pi B_{0, N} (z)^{2}}.
\end{equation*}
\end{corollary}

We shall use the formula for $h_{N, \mathbf{K}} (z)$ in Theorem \ref{theorem-1} for the special choices of $f_{j} (z)$ to study the limiting behaviour of $h_{N, \mathbf{K}} (z)$ as $N$ tends to infinity. This shall demonstrate how the zeros of the random equation $S_{N} (z) = \mathbf{K}$ are clustered in the limit. Finally, using the Christoffel--Darboux formula, we derive the density function for the complex zeros of orthogonal polynomials, with the orthogonality relation being satisfied on the real number line and the unit circle.

\medskip

\section{The evaluation of the density function}

We shall begin the proof of Theorem \ref{theorem-1} by letting $X_{1, N}$ and $X_{2, N}$ be the real and imaginary parts of $S_{N} (z)$, respectively. For convenience in computation, we shall write
\begin{equation*}
f_{j} (z)
 = u_{j} (x, y) + i v_{j} (x, y),
\end{equation*}
where $u_{j} (x, y)$ and $v_{j} (x, y)$ are real-valued functions of $(x, y) \in \mathds{R}^{2}$. We have
\begin{equation*}
S_{N} (z)
 = X_{1, N} + i X_{2, N},
\end{equation*}
where
\begin{equation*}
X_{1, N}
 = \sum_{j = 0}^{N} (a_{j} u_{j} - b_{j} v_{j})
\end{equation*}
and
\begin{equation*}
X_{2, N}
 = \sum_{j = 0}^{N} (a_{j} v_{j} + b_{j} u_{j}).
\end{equation*}
In our application of Adler's theorem, we need to find all real and complex zeros of $S_{N} (z)$. They are the zeros of the random equations $X_{1, N} = K_{1}$ and $X_{2, N} = K_{2}$ for $(x, y) \in \mathds{R}^{2}$.

For the sake of brevity, we let $\mathbf{X}_{N}$ be the two-dimensional random field of the real and imaginary parts of $S_{N} (z)$ defined by the column vector $\mathbf{X}_{N} = (X_{1, N}, X_{2, N})^{\prime}$, and we denote the Jacobian matrix of the transformation $(x, y) \longrightarrow (X_{1, N}, X_{2, N})$ by the matrix $\nabla \mathbf{X}_{N}$ of the first-order partial derivatives of $\mathbf{X}_{N}$ with respect to $x$ and $y$, namely,
\begin{equation*}
\nabla \mathbf{X}_{N}
 = \frac{\partial (X_{1, N}, X_{2, N})}{\partial (x, y)}
 =
\begin{pmatrix}
\dfrac{\partial X_{1, N}}{\partial x} & \dfrac{\partial X_{2, N}}{\partial x} \\ \addlinespace \addlinespace \addlinespace
\dfrac{\partial X_{1, N}}{\partial y} & \dfrac{\partial X_{2, N}}{\partial y}
\end{pmatrix}.
\end{equation*}
We find the determinant of $\nabla \mathbf{X}_{N}$. It will be convenient to first use the Cauchy--Riemann equations to rewrite the expressions for $\partial X_{1, N} / \partial y$ and $\partial X_{2, N} / \partial y$. In order to obtain the conditional expectation of $\abs{\mbox{det} \nabla \mathbf{X}_{N}}$ on the extreme right side of \eqref{equation-3}, we separate the diagonal terms from the cross terms in the random determinant $\mbox{det} \nabla \mathbf{X}_{N}$. It is an easy computation, and with a little algebra we find that
\begin{align} \label{equation-2}
\mbox{det} \nabla \mathbf{X}_{N}
 &= 
\begin{vmatrix}[l]
\, \displaystyle \sum_{j = 0}^{N} \left(a_{j} \frac{\partial u_{j}}{\partial x} - b_{j} \frac{\partial v_{j}}{\partial x}\right) & \displaystyle \sum_{j = 0}^{N} \left(a_{j} \frac{\partial v_{j}}{\partial x} + b_{j} \frac{\partial u_{j}}{\partial x}\right) \, \\ \addlinespace \addlinespace \nonumber
\, \displaystyle \sum_{j = 0}^{N} \left(-a_{j} \frac{\partial v_{j}}{\partial x} - b_{j} \frac{\partial u_{j}}{\partial x}\right) & \displaystyle \sum_{j = 0}^{N} \left(a_{j} \frac{\partial u_{j}}{\partial x} - b_{j} \frac{v_{j}}{\partial x}\right) \,
\end{vmatrix} \\ \nonumber
 &= \sum_{j = 0}^{N} \sum_{k = 0}^{N} \left\{(a_{j} a_{k} + b_{j} b_{k}) \left(\frac{\partial u_{j}}{\partial x} \frac{\partial u_{k}}{\partial x} + \frac{\partial v_{j}}{\partial x} \frac{\partial v_{k}}{\partial x}\right)\right. \nonumber \\ &\hspace{95pt} \left. + (a_{j} b_{k} - b_{j} a_{k}) \left(\frac{\partial v_{j}}{\partial x} \frac{\partial u_{k}}{\partial x} - \frac{\partial u_{j}}{\partial x} \frac{\partial v_{k}}{\partial x}\right)\right\} \nonumber \\
 &= \sum_{j = 0}^{N} (a_{j}^{2} + b_{j}^{2}) \left\{\left(\frac{\partial u_{j}}{\partial x}\right)^{2} + \left(\frac{\partial v_{j}}{\partial x}\right)^{2}\right\} \nonumber \\ &\quad + \sum_{j = 0}^{N} \sum_{\substack{k = 0 \\ k \neq j}}^{N} \left\{(a_{j} a_{k} + b_{j} b_{k}) \left(\frac{\partial u_{j}}{\partial x} \frac{\partial u_{k}}{\partial x} + \frac{\partial v_{j}}{\partial x} \frac{\partial v_{k}}{\partial x}\right)\right. \nonumber \\ &\hspace{95pt} \left. + (a_{j} b_{k} - b_{j} a_{k}) \left(\frac{\partial v_{j}}{\partial x} \frac{\partial u_{k}}{\partial x} - \frac{\partial u_{j}}{\partial x} \frac{\partial v_{k}}{\partial x}\right)\right\}.
\end{align}
Thus, the evaluation of $h_{N, \mathbf{K}} (z)$ leads to the computation of the expected value of a quadratic form $\mbox{det} \nabla \mathbf{X}_{N}$ of independent, identically distributed Gaussian random variables, conditioned on two linear combinations.

We shall assume that there are no points on the boundary $\partial \Phi$ for which $\mathbf{X}_{N} = \mathbf{0}$, and furthermore that the compact subset $\Phi \subset \mathds{C}$ does not contain any points satisfying $\mathbf{X}_{N} = \mathbf{0}$ and $\nabla \mathbf{X}_{N} = \mathbf{0}$ at the same time. Then by Theorem 5.1.1 and its Corollary in Adler's classical book \cite[pages 95--97]{Adler1981} (see, also, the papers by Aza\"{i}s and Wschebor \cite{AzaisWschebor2005} and Ibragimov and Zeitouni \cite{IbragimovZeitouni1997}), the density function $h_{N, \mathbf{K}} (z)$ for multivariate Gaussian coefficients given by \eqref{equation-1} can be expressed through a conditioned expected value given by
\begin{align} \label{equation-3}
h_{N, \mathbf{K}} (z)
 &= \mathscr{E} (\abs{\mbox{det} \nabla \mathbf{X}_{N}} \mid X_{1, N} = K_{1}, X_{2, N} = K_{2}) \, p_{x, y} (K_{1}, K_{2}) \nonumber \\
 &= \mathscr{E} (\abs{\mbox{det} \nabla \mathbf{X}_{N}} \mid \mathbf{X}_{N} = \mathbf{K}) \, p_{x, y} (\mathbf{K}^{\prime}),
\end{align}
where $p_{x, y} (\mathbf{K}^{\prime})$ is the two-dimensional joint density function of the random vector $\mathbf{X}_{N}$. Since $\mbox{det} \nabla \mathbf{X}_{N}$ is always nonnegative, let us eliminate the absolute value sign from future occurrences of the extreme right side of \eqref{equation-3} in the evaluation of $h_{N, \mathbf{K}} (z)$.

In the following, we obtain the vectors of conditional expectations, variances, and covariance matrices of the multivariate random vectors $\mathbf{a}_{N} = (a_{0}, \ldots, a_{N})^{\prime}$ and $\mathbf{b}_{N} = (b_{0}, \ldots, b_{N})^{\prime}$. From standard methods in multivariate analysis (see the classical books by Anderson \cite{Anderson1984} and Tong \cite{Tong1990}), based on the assumption that all the scalar random variables involved are independent and normally distributed, we define
\begin{equation} \label{equation-4}
\mbox{Cov}
(\mathbf{a}_N, \mathbf{b}_{N} \mid \mathbf{X}_{N} = \mathbf{K})
 = 
\begin{pmatrix}
\, \mathbf{\Sigma}_{\mathbf{a}_{N} \mathbf{a}_{N}, \mathbf{X}_{N}} & \mathbf{\Sigma}_{\mathbf{a}_{N} \mathbf{b}_{N}, \mathbf{X}_{N}} \, \\ \addlinespace
\, \mathbf{\Sigma}_{\mathbf{b}_{N} \mathbf{a}_{N}, \mathbf{X}_{N}} & \mathbf{\Sigma}_{\mathbf{b}_{N} \mathbf{b}_{N}, \mathbf{X}_{N}} \,
\end{pmatrix}
\end{equation}
and
\begin{equation} \label{equation-5}
\mathscr{E} (\mathbf{a}_{N} \mid \mathbf{X}_{N} = \mathbf{K})
 = \mathscr{E} \mathbf{a}_{N} + \mathbf{\Sigma}_{\mathbf{a}_{N} \mathbf{X}_{N}} \mathbf{\Sigma}_{\mathbf{X}_{N} \mathbf{X}_{N}}^{-1} (\mathbf{K} - \mathscr{E} \mathbf{X}_{N})^{\prime},
\end{equation}
where
\begin{equation} \label{equation-6}
\mathbf{\Sigma}_{\mathbf{a}_{N} \mathbf{b}_{N}, \mathbf{X}_{N}}
 = \mathbf{\Sigma}_{\mathbf{a}_{N} \mathbf{b}_{N}} - \mathbf{\Sigma}_{\mathbf{a}_{N} \mathbf{X}_{N}} \mathbf{\Sigma}_{\mathbf{X}_{N} \mathbf{X}_{N}}^{-1} \mathbf{\Sigma}_{\mathbf{X}_{N} \mathbf{b}_{N}}
\end{equation}
and
\begin{equation} \label{equation-7}
\mathbf{\Sigma}_{\mathbf{a}_{N} \mathbf{b}_{N}}
 = \mathscr{E} (\mathbf{a}_{N} - \mathscr{E} \mathbf{a}_{N}) (\mathbf{b}_{N} - \mathscr{E} \mathbf{b}_{N})^{\prime},
\end{equation}
which is a generalized covariance matrix of the vectors $\mathbf{a}_{N}$ and $\mathbf{b}_{N}$. Whereas the distribution of the $a_{j}$ and $b_{j}$ is central, we have $\mathscr{E} \mathbf{a}_{N} = 0$, $\mathscr{E} \mathbf{b}_{N} = 0$, and $\mathscr{E} \mathbf{X}_{N} = 0$.

From \eqref{equation-7} and the assumption of the theorem
\begin{equation} \label{equation-8}
\mathbf{\Sigma}_{\mathbf{a}_{N} \mathbf{a}_{N}}
 = \mathscr{E} (\mathbf{a}_{N} - \mathscr{E} \mathbf{a}_{N}) (\mathbf{a}_{N} - \mathscr{E} \mathbf{a}_{N})^{\prime}
 = \mathscr{E} \mathbf{a}_{N} \mathbf{a}_{N}^{\prime}
 = \mathbf{I}_{N},
\end{equation}
since the $a_{j}$ are distributed according to an $\mathscr{N} (0, 1)$ distribution and
\begin{equation} \label{equation-9}
\mathscr{E} a_{j} a_{k}
 = \left\{ \begin{array}{ll}
   1 & \mbox{if $j = k$,} \\ \addlinespace
   0 & \mbox{if $j \neq k$.}
\end{array}
\right.
\end{equation}
In a similar fashion, from \eqref{equation-7}
\begin{equation} \label{equation-10}
\mathbf{\Sigma}_{\mathbf{b}_{N} \mathbf{b}_{N}}
 = \mathscr{E} (\mathbf{b}_{N} - \mathscr{E} \mathbf{b}_{N}) (\mathbf{b}_{N} - \mathscr{E} \mathbf{b}_{N})^{\prime}
 = \mathscr{E} \mathbf{b}_{N} \mathbf{b}_{N}^{\prime}
 = \mathbf{I}_{N}.
\end{equation}
Since
\begin{equation} \label{equation-11}
\mathscr{E} a_{j} b_{k}
 = 0,
\end{equation}
we have
\begin{equation} \label{equation-12}
\mathbf{\Sigma}_{\mathbf{a}_{N} \mathbf{b}_{N}}
 = \mathscr{E} (\mathbf{a}_{N} - \mathscr{E} \mathbf{a}_{N}) (\mathbf{b}_{N} - \mathscr{E} \mathbf{b}_{N})^{\prime}
 = \mathscr{E} \mathbf{a}_{N} \mathbf{b}_{N}^{\prime} = \mathbf{0}.
\end{equation}
It follows that
\begin{equation} \label{equation-13}
\mathbf{\Sigma}_{\mathbf{b}_{N} \mathbf{a}_{N}}
 = \mathbf{\Sigma}_{\mathbf{a}_{N} \mathbf{b}_{N}}^{\prime} = \mathbf{0}.
\end{equation}

Now, from \eqref{equation-9} and \eqref{equation-11}
\begin{equation} \label{equation-14}
\mathscr{E} a_{j} X_{1, N}
 = \sum_{k = 0}^{N} (\mathscr{E} (a_{j} a_{k}) u_{k} - \mathscr{E} (a_{j} b_{k}) v_{k})
 = u_{j}
\end{equation} 
and
\begin{equation} \label{equation-15}
\mathscr{E} a_{j} X_{2, N}
 = \sum_{k = 0}^{N} (\mathscr{E} (a_{j} a_{k}) v_{k} + \mathscr{E} (a_{j} b_{k}) u_{k})
 = v_{j}.
\end{equation}
If we apply \eqref{equation-14} and \eqref{equation-15} to \eqref{equation-7}, we obtain
\begin{align} \label{equation-16}
\mathbf{\Sigma}_{\mathbf{a}_{N} \mathbf{X}_{N}}
 &= \mathscr{E} (\mathbf{a}_{N} - \mathscr{E} \mathbf{a}_{N}) (\mathbf{X}_{N} - \mathscr{E} \mathbf{X}_{N})^{\prime}
  = \mathscr{E} \mathbf{a}_{N} \mathbf{X}_{N}^{\prime} \nonumber \\
 &=
\begin{pmatrix}
 \mathscr{E} a_{0} X_{1, N}	& \mathscr{E} a_{0} X_{2, N} \\ \addlinespace
 \mathscr{E} a_{1} X_{1, N}	& \mathscr{E} a_{1} X_{2, N} \\ \addlinespace
 \hdotsfor{2}				 \\ \addlinespace
\mathscr{E} a_{N} X_{1, N}	& \mathscr{E} a_{N} X_{2, N}
\end{pmatrix}
=
\begin{pmatrix}
 u_{0}	& v_{0} \\ \addlinespace
 u_{1}	& v_{1} \\ \addlinespace
 \hdotsfor{2} \\ \addlinespace
 u_{N}	& v_{N}
\end{pmatrix}.
\end{align}
Thus,
\begin{equation} \label{equation-17}
\mathbf{\Sigma}_{\mathbf{X}_{N} \mathbf{a}_{N}}
 = \mathbf{\Sigma}_{\mathbf{a}_{N} \mathbf{X}_{N}}^{\prime}
=
\begin{pmatrix}
 u_{0}	& u_{1} 	& \hdotsfor{1}	& u_{N} \\ \addlinespace
 v_{0}	& v_{1} 	& \dots	& v_{N}
\end{pmatrix}.
\end{equation}
Proceeding as above, using \eqref{equation-9} and \eqref{equation-11} with the obvious substitutions, we obtain
\begin{equation} \label{equation-18}
\mathscr{E} b_{j} X_{1, N}
 = \sum_{k = 0}^{N} (\mathscr{E} (b_{j} a_{k}) u_{k} - \mathscr{E} (b_{j} b_{k}) v_{k})
 = -v_{j}
\end{equation}
and
\begin{equation} \label{equation-19}
\mathscr{E} b_{j} X_{2, N}
 = \sum_{k = 0}^{N} (\mathscr{E} (b_{j} a_{k}) v_{k} + \mathscr{E} (b_{j} b_{k}) u_{k})
 = u_{j}.
\end{equation}
Then applying \eqref{equation-18} and \eqref{equation-19} to \eqref{equation-7}, we get
\begin{align} \label{equation-20}
\mathbf{\Sigma}_{\mathbf{b}_{N} \mathbf{X}_{N}}
 &= \mathscr{E} (\mathbf{b}_{N} - \mathscr{E} \mathbf{b}_{N}) (\mathbf{X}_{N} - \mathscr{E} \mathbf{X}_{N})^{\prime} 
  = \mathscr{E} \mathbf{b}_{N} \mathbf{X}_{N}^{\prime} \nonumber \\
 &=
\begin{pmatrix}
 \mathscr{E} b_{0} X_{1, N}	& \mathscr{E} b_{0} X_{2, N} \\ \addlinespace
 \mathscr{E} b_{1} X_{1, N}	& \mathscr{E} b_{1} X_{2, N} \\ \addlinespace
 \hdotsfor{2}				 \\ \addlinespace
 \mathscr{E} b_{N} X_{1, N}	& \mathscr{E} b_{N} X_{2, N}
\end{pmatrix}
=
\begin{pmatrix}
 -v_{0}	& u_{0} \\ \addlinespace
 -v_{1}	& u_{1} \\ \addlinespace
 \hdotsfor{2} \\ \addlinespace
 -v_{N}	& u_{N}
\end{pmatrix}.
\end{align}
Hence,
\begin{equation} \label{equation-21}
\mathbf{\Sigma}_{\mathbf{X}_{N} \mathbf{b}_{N}}
 = \mathbf{\Sigma}_{\mathbf{b}_{N} \mathbf{X}_{N}}^{\prime}
=
\begin{pmatrix}
 -v_{0}	& -v_{1}	& \dots	& -v_{N} \\ \addlinespace
 u_{0}	& u_{1}	& \dots	& u_{N}
\end{pmatrix}.
\end{equation}

Again, we use \eqref{equation-7} to obtain
\begin{align} \label{equation-22}
\mathbf{\Sigma}_{\mathbf{X}_{N} \mathbf{X}_{N}}
 &= \mathscr{E} (\mathbf{X}_{N} - \mathscr{E} \mathbf{X}_{N}) (\mathbf{X}_{N} - \mathscr{E} \mathbf{X}_{N})^{\prime} \nonumber \\
 &=
\begin{pmatrix}
\mathscr{E} X_{1, N} X_{1, N}	& \mathscr{E} X_{1, N} X_{2, N} \\ \addlinespace
\mathscr{E} X_{2, N} X_{1, N}	& \mathscr{E} X_{2, N} X_{2, N}
\end{pmatrix}.
\end{align}
We compute that
\begin{align} \label{equation-23}
\mathscr{E} X_{1, N} X_{1, N}
 &= \mathscr{E} \left(\sum_{j = 0}^{N} \sum_{k = 0}^{N} (a_{j} u_{j} - b_{j} v_{j}) (a_{k} u_{k} - b_{k} v_{k})\right) \nonumber \\
 &= \sum_{j = 0}^{N} \sum_{k = 0}^{N} (\mathscr{E} (a_{j} a_{k}) u_{j} u_{k} + \mathscr{E} (b_{j} b_{k}) v_{j} v_{k} - \mathscr{E} (a_{j} b_{k}) u_{j} v_{k} - \mathscr{E} (b_{j} a_{k}) v_{j} u_{k}) \nonumber \\
 &= \sum_{j = 0}^{N} (u_{j}^{2} + v_{j}^{2})
 = \sum_{j = 0}^{N} \abs{f_{j} (z)}^{2}.
\end{align}
We get, similarly to \eqref{equation-23},
\begin{align} \label{equation-24}
\mathscr{E} X_{2, N} X_{2, N}
 &=\mathscr{E} \left(\sum_{j = 0}^{N} \sum_{k = 0}^{N} (a_{j} v_{j} + b_{j}u_{j}) (a_{k} v_{k} + b_{k} v_{k})\right) \nonumber \\
 &= \sum_{j = 0}^{N} \sum_{k = 0}^{N} (\mathscr{E} (a_{j} a_{k}) v_{j} v_{k} + \mathscr{E} (b_{j} b_{k}) u_{j} u_{k} + \mathscr{E} (a_{j} b_{k}) v_{j} u_{k} + \mathscr{E} (b_{j} a_{k}) u_{j} v_{k}) \nonumber \\
 &= \sum_{j = 0}^{N} (u_{j}^{2} + v_{j}^{2})
 = \sum_{j = 0}^{N} \abs{f_{j} (z)}^{2}.
\end{align}
Furthermore, we have
\begin{align} \label{equation-25}
\mathscr{E} X_{2, N} X_{1, N}
 &= \mathscr{E} X_{1, N} X_{2, N}
 = \mathscr{E} \left(\sum_{j = 0}^{N} \sum_{k = 0}^{N} (a_{j} u_{j} - b_{j} v_{j}) (a_{k} v_{k} + b_{k} u_{k})\right) \nonumber \\
 &= \sum_{j = 0}^{N} \sum_{k = 0}^{N} \left(\mathscr{E} (a_{j} a_{k}) u_{j} v_{k} - \mathscr{E} (b_{j} b_{k}) v_{j} u_{k} + \mathscr{E} (a_{j} b_{k}) u_{j} u_{k} - \mathscr{E} (b_{j} a_{k}) v_{j} v_{k}\right) \nonumber \\
 &= \sum_{j = 0}^{N} (u_{j} v_{j} - v_{j} u_{j})
 = 0.
\end{align}
Then combining \eqref{equation-23}--\eqref{equation-25} in \eqref{equation-22}, we obtain
\begin{equation*}
\mathbf{\Sigma}_{\mathbf{X}_{N} \mathbf{X}_{N}}
= \sum_{j = 0}^{N} \abs{f_{j} (z)}^{2} \mathbf{I}_{2}.
\end{equation*}
We note that the existence of $h_{N, \mathbf{K}} (z)$ depends on the evaluation of the covariance matrix $\mbox{Cov}
(\mathbf{a}_N, \mathbf{b}_{N} \mid \mathbf{X}_{N} = \mathbf{K})$, which in turn depends on the existence of the inverse matrix $\mathbf{\Sigma}_{\mathbf{X}_{N} \mathbf{X}_{N}}^{-1}$. This is guaranteed, since
\begin{equation*}
\mbox{det} \abs{\mathbf{\Sigma}_{\mathbf{X}_{N} \mathbf{X}_{N}}}
 = \left(\sum_{j = 0}^{N} \abs{f_{j} (z)}^{2}\right)^{2}.
\end{equation*}
Thus,
\begin{equation} \label{equation-26}
\mathbf{\Sigma}_{\mathbf{X}_{N} \mathbf{X}_{N}}^{-1}
 = \frac{\mathbf{I}_{2}}{\sqrt{\mbox{det} \abs{\mathbf{\Sigma}_{\mathbf{X}_{N} \mathbf{X}_{N}}}}}.
\end{equation}

Moving now to the components of the covariance matrix given by \eqref{equation-4}, we obtain the results for the $j$th row and the $k$th column. Let $\delta_{j k}$ denote the Kronecker delta, that is,
\begin{equation*}
\delta_{j k}
 = \left\{ \begin{array}{ll}
   1 & \mbox{if $j = k$,} \\
   0 & \mbox{if $j \neq k$.}
\end{array}
\right.
\end{equation*}
From \eqref{equation-8}, \eqref{equation-16}, \eqref{equation-17}, and \eqref{equation-26}
\begin{align} \label{equation-27}
(\mathbf{\Sigma}_{\mathbf{a}_{N} \mathbf{a}_{N}, \mathbf{X}_{N}})_{\substack{1 \leq j \leq N \\ 1 \leq k \leq N}}
 &= (\mathbf{\Sigma}_{\mathbf{a}_{N} \mathbf{a}_{N}})_{\substack{1 \leq j \leq N \\ 1 \leq k \leq N}} - (\mathbf{\Sigma}_{\mathbf{a}_{N} \mathbf{X}_{N}} \mathbf{\Sigma}_{\mathbf{X}_{N} \mathbf{X}_{N}}^{-1} \mathbf{\Sigma}_{\mathbf{X}_{N} \mathbf{a}_{N}})_{\substack{1 \leq j \leq N \\ 1 \leq k \leq N}} \nonumber \\
 &= \delta_{j k} - \frac{u_{j} u_{k} + v_{j} v_{k}}{\sqrt{\mbox{det} \abs{\mathbf{\Sigma}_{\mathbf{X}_{N} \mathbf{X}_{N}}}}}.
\end{align}
From \eqref{equation-10}, \eqref{equation-20}, \eqref{equation-21}, and \eqref{equation-26}
\begin{align} \label{equation-28}
(\mathbf{\Sigma}_{\mathbf{b}_{N} \mathbf{b}_{N}, \mathbf{X}_{N}})_{\substack{1 \leq j \leq N \\ 1 \leq k \leq N}}
 &= (\mathbf{\Sigma}_{\mathbf{b}_{N} \mathbf{b}_{N}})_{\substack{1 \leq j \leq N \\ 1 \leq k \leq N}} - (\mathbf{\Sigma}_{\mathbf{b}_{N} \mathbf{X}_{N}} \mathbf{\Sigma}_{\mathbf{X}_{N} \mathbf{X}_{N}}^{-1} \mathbf{\Sigma}_{\mathbf{X}_{N} \mathbf{b}_{N}})_{\substack{1 \leq j \leq N \\ 1 \leq k \leq N}} \nonumber \\
 &= \delta_{j k} - \frac{u_{j} u_{k} + v_{j} v_{k}}{\sqrt{\mbox{det} \abs{\mathbf{\Sigma}_{\mathbf{X}_{N} \mathbf{X}_{N}}}}}.
\end{align}
From \eqref{equation-12}, \eqref{equation-16}, \eqref{equation-21}, and \eqref{equation-26}
\begin{align} \label{equation-29}
(\mathbf{\Sigma}_{\mathbf{a}_{N} \mathbf{b}_{N}, \mathbf{X}_{N}})_{\substack{1 \leq j \leq N \\ 1 \leq k \leq N}}
 &= (\mathbf{\Sigma}_{\mathbf{a}_{N} \mathbf{b}_{N}})_{\substack{1 \leq j \leq N \\ 1 \leq k \leq N}} - (\mathbf{\Sigma}_{\mathbf{a}_{N} \mathbf{X}_{N}} \mathbf{\Sigma}_{\mathbf{X}_{N} \mathbf{X}_{N}}^{-1} \mathbf{\Sigma}_{\mathbf{X}_{N} \mathbf{b}_{N}})_{\substack{1 \leq j \leq N \\ 1 \leq k \leq N}} \nonumber \\
 &= \frac{u_{j} v_{k} - v_{j} u_{k}}{\sqrt{\mbox{det} \abs{\mathbf{\Sigma}_{\mathbf{X}_{N} \mathbf{X}_{N}}}}}.
\end{align}
From \eqref{equation-13}, \eqref{equation-17}, \eqref{equation-20}, and \eqref{equation-26}
\begin{align} \label{equation-30}
(\mathbf{\Sigma}_{\mathbf{b}_{N} \mathbf{a}_{N}, \mathbf{X}_{N}})_{\substack{1 \leq j \leq N \\ 1 \leq k \leq N}}
 &= (\mathbf{\Sigma}_{\mathbf{b}_{N} \mathbf{a}_{N}})_{\substack{1 \leq j \leq N \\ 1 \leq k \leq N}} - (\mathbf{\Sigma}_{\mathbf{b}_{N} \mathbf{X}_{N}} \mathbf{\Sigma}_{\mathbf{X}_{N} \mathbf{X}_{N}}^{-1} \mathbf{\Sigma}_{\mathbf{X}_{N} \mathbf{a}_{N}})_{\substack{1 \leq j \leq N \\ 1 \leq k \leq N}} \nonumber \\
 &= -\frac{u_{j} v_{k} - v_{j} u_{k}}{\sqrt{\mbox{det} \abs{\mathbf{\Sigma}_{\mathbf{X}_{N} \mathbf{X}_{N}}}}}.
\end{align}

We proceed to finding the necessary conditional expectations for computing $\mathscr{E} (\mbox{det} \nabla \mathbf{X}_{N} \mid \mathbf{X}_{N} = \mathbf{K})$. The conditional expectation of $\mathbf{a}$ and $\mathbf{b}$ is easily derived, respectively, from \eqref{equation-16} and \eqref{equation-26} as
\begin{align} \label{equation-31}
(\mathscr{E} (a_{j} \mid \mathbf{X}_{N} = \mathbf{K}))_{1 \leq j \leq N}
 &= (\mathbf{\Sigma}_{\mathbf{a}_{N} \mathbf{X}_{N}} \mathbf{\Sigma}_{\mathbf{X}_{N} \mathbf{X}_{N}}^{-1} \mathbf{K})_{1 \leq j \leq N} \nonumber \\
 &= \frac{K_{1} u_{j} + K_{2} v_{j}}{\sqrt{\mbox{det} \abs{\mathbf{\Sigma}_{\mathbf{X}_{N} \mathbf{X}_{N}}}}}
\end{align}
and from \eqref{equation-20} and \eqref{equation-26} as
\begin{align} \label{equation-32}
(\mathscr{E} (b_{j} \mid \mathbf{X}_{N} = \mathbf{K}))_{1 \leq j \leq N}
 &= (\mathbf{\Sigma}_{\mathbf{b}_{N} \mathbf{X}_{N}} \mathbf{\Sigma}_{\mathbf{X}_{N} \mathbf{X}_{N}}^{-1} \mathbf{K})_{1 \leq j \leq N} \nonumber \\
 &= \frac{K_{2} u_{j} - K_{1} v_{j}}{\sqrt{\mbox{det} \abs{\mathbf{\Sigma}_{\mathbf{X}_{N} \mathbf{X}_{N}}}}}.
\end{align}
We derive from \eqref{equation-31}
\begin{align} \label{equation-33}
(\mathscr{E} (a_{j}^{2} \mid \mathbf{X}_{N} = \mathbf{K}))_{1 \leq j \leq N}
 &= (\mathscr{E} (a_{j} \mid \mathbf{X}_{N} = \mathbf{K}))_{1 \leq j \leq N}^{2} + (\mbox{Var} (a_{j} \mid \mathbf{X}_{N} = \mathbf{K}))_{1 \leq j \leq N} \nonumber \\
 &= \frac{K_{1}^{2} u_{j}^{2} + K_{2}^{2} v_{j}^{2} + 2 K_{1} K_{2} v_{j} v_{j}}{\mbox{det} \abs{\mathbf{\Sigma}_{\mathbf{X}_{N} \mathbf{X}_{N}}}} + 1 \nonumber \\ &\quad - \frac{u_{j}^{2} + v_{j}^{2}}{\sqrt{\mbox{det} \abs{\mathbf{\Sigma}_{\mathbf{X}_{N} \mathbf{X}_{N}}}}}
\end{align}
and from \eqref{equation-32}
\begin{align} \label{equation-34}
(\mathscr{E} (b_{j}^{2} \mid \mathbf{X}_{N} = \mathbf{K}))_{1 \leq j \leq N}
 &= (\mathscr{E} (b_{j} \mid \mathbf{X}_{N} = \mathbf{K}))_{1 \leq j \leq N}^{2} + (\mbox{Var} (b_{j} \mid \mathbf{X}_{N} = \mathbf{K}))_{1 \leq j \leq N} \nonumber \\
 &= \frac{K_{2}^{2} u_{j}^{2} + K_{1}^{2} v_{j}^{2} - 2 K_{1} K_{2} v_{j} v_{j}}{\mbox{det} \abs{\mathbf{\Sigma}_{\mathbf{X}_{N} \mathbf{X}_{N}}}} + 1 \nonumber \\ &\quad - \frac{u_{j}^{2} + v_{j}^{2}}{\sqrt{\mbox{det} \abs{\mathbf{\Sigma}_{\mathbf{X}_{N} \mathbf{X}_{N}}}}}.
\end{align}
By virtue of \eqref{equation-33} and \eqref{equation-34}
\begin{align} \label{equation-35}
(\mathscr{E} (a_{j}^{2} + b_{j}^{2} \mid \mathbf{X}_{N} = \mathbf{K}))_{1 \leq j \leq N}
 &= \frac{(K_{1}^{2} + K_{2}^{2}) (u_{j}^{2} + v_{j}^{2})}{\mbox{det} \abs{\mathbf{\Sigma}_{\mathbf{X}_{N} \mathbf{X}_{N}}}} + 2 \nonumber \\ &\quad - \frac{2 (u_{j}^{2} + v_{j}^{2})}{\sqrt{\mbox{det} \abs{\mathbf{\Sigma}_{\mathbf{X}_{N} \mathbf{X}_{N}}}}}.
\end{align}
Next, using \eqref{equation-4} and \eqref{equation-27}--\eqref{equation-31}, we get
\begin{align} \label{equation-36}
(\mathscr{E} (a_{j} a_{k} \mid \mathbf{X}_{N} = \mathbf{K}))_{\substack{1 \leq j \leq N \\ 1 \leq k \leq N}}
 &= (\mathscr{E} (a_{j} \mid \mathbf{X}_{N} = \mathbf{K}))_{1 \leq j \leq N} \, (\mathscr{E} (a_{k} \mid \mathbf{X}_{N} = \mathbf{K}))_{1 \leq k \leq N} \nonumber \\ &\quad + (\mbox{Cov} (a_{j}, a_{k} \mid \mathbf{X}_{N} = \mathbf{K}))_{\substack{1 \leq j \leq N \\ 1 \leq k \leq N}} \nonumber \\
 &= \frac{K_{1}^{2} u_{j} u_{k} + K_{2}^{2} v_{j} v_{k} + K_{1} K_{2} (u_{j} v_{k} + v_{j} u_{k})}{\mbox{det} \abs{\mathbf{\Sigma}_{\mathbf{X}_{N} \mathbf{X}_{N}}}} \nonumber \\ &\quad - \frac{u_{j} u_{k} + v_{j} v_{k}}{\sqrt{\mbox{det} \abs{\mathbf{\Sigma}_{\mathbf{X}_{N} \mathbf{X}_{N}}}}}.
\end{align}
Using \eqref{equation-4}, \eqref{equation-27}--\eqref{equation-30}, and \eqref{equation-32}, we find
\begin{align} \label{equation-37}
(\mathscr{E} (b_{j} b_{k} \mid \mathbf{X}_{N} = \mathbf{K}))_{\substack{1 \leq j \leq N \\ 1 \leq k \leq N}}
 &= (\mathscr{E} (b_{j} \mid \mathbf{X}_{N} = \mathbf{K}))_{1 \leq j \leq N} \, (\mathscr{E} (b_{k} \mid \mathbf{X}_{N} = \mathbf{K}))_{1 \leq k \leq N} \nonumber \\ &\quad + (\mbox{Cov} (b_{j}, b_{k} \mid \mathbf{X}_{N} = \mathbf{K}))_{\substack{1 \leq j \leq N \\ 1 \leq k \leq N}} \nonumber \\
 &= \frac{K_{2}^{2} u_{j} u_{k} + K_{1}^{2} v_{j} v_{k} - K_{1} K_{2} (u_{j} v_{k} + v_{j} u_{k})}{\mbox{det} \abs{\mathbf{\Sigma}_{\mathbf{X}_{N} \mathbf{X}_{N}}}} \nonumber \\ &\quad - \frac{u_{j} u_{k} + v_{j} v_{k}}{\sqrt{\mbox{det} \abs{\mathbf{\Sigma}_{\mathbf{X}_{N} \mathbf{X}_{N}}}}}.
\end{align}
By virtue of \eqref{equation-36} and \eqref{equation-37}
\begin{align} \label{equation-38}
(\mathscr{E} (a_{j} a_{k} + b_{j} b_{k} \mid \mathbf{X}_{N} = \mathbf{K}))_{\substack{1 \leq j \leq N \\ 1 \leq k \leq N}}
 &= \frac{(K_{1}^{2} + K_{2}^{2}) (u_{j} u_{k} + v_{j} v_{k})}{\mbox{det} \abs{\mathbf{\Sigma}_{\mathbf{X}_{N} \mathbf{X}_{N}}}} \nonumber \\ &\quad - \frac{2 (u_{j} u_{k} + v_{j} v_{k})}{\sqrt{\mbox{det} \abs{\mathbf{\Sigma}_{\mathbf{X}_{N} \mathbf{X}_{N}}}}}.
\end{align}

Next, we derive from \eqref{equation-4} and \eqref{equation-27}--\eqref{equation-32}
\begin{align} \label{equation-39}
(\mathscr{E} (a_{j} b_{k} \mid \mathbf{X}_{N} = \mathbf{K}))_{\substack{1 \leq j \leq N \\ 1 \leq k \leq N}}
 &= (\mathscr{E} (a_{j} \mid \mathbf{X}_{N} = \mathbf{K}))_{1 \leq j \leq N} \, (\mathscr{E} (b_{k} \mid \mathbf{X}_{N} = \mathbf{K}))_{1 \leq k \leq N} \nonumber \\ &\quad + (\mbox{Cov} (a_{j}, b_{k} \mid \mathbf{X}_{N} = \mathbf{K}))_{\substack{1 \leq j \leq N \\ 1 \leq k \leq N}} \nonumber \\
 &= \frac{K_{1} K_{2} (u_{j} u_{k} - v_{j} v_{k}) - K_{1}^{2} u_{j} v_{k} + K_{2}^{2} v_{j} u_{k}}{\mbox{det} \abs{\mathbf{\Sigma}_{\mathbf{X}_{N} \mathbf{X}_{N}}}} \nonumber \\ &\quad + \frac{u_{j} v_{k} - v_{j} u_{k}}{\sqrt{\mbox{det} \abs{\mathbf{\Sigma}_{\mathbf{X}_{N} \mathbf{X}_{N}}}}}
\end{align}
and
\begin{align} \label{equation-40}
(\mathscr{E} (b_{j} a_{k} \mid \mathbf{X}_{N} = \mathbf{K}))_{\substack{1 \leq j \leq N \\ 1 \leq k \leq N}}
 &= (\mathscr{E} (b_{j} \mid \mathbf{X}_{N} = \mathbf{K}))_{1 \leq j \leq N} \, (\mathscr{E} (a_{k} \mid \mathbf{X}_{N} = \mathbf{K}))_{1 \leq k \leq N} \nonumber \\ &\quad + (\mbox{Cov} (b_{j}, a_{k} \mid \mathbf{X}_{N} = \mathbf{K}))_{\substack{1 \leq j \leq N \\ 1 \leq k \leq N}} \nonumber \\
 &= \frac{K_{1} K_{2} (u_{j} u_{k} - v_{j} v_{k}) + K_{2}^{2} u_{j} v_{k} - K_{1}^{2} v_{j} u_{k}}{\mbox{det} \abs{\mathbf{\Sigma}_{\mathbf{X}_{N} \mathbf{X}_{N}}}} \nonumber \\ &\quad - \frac{u_{j} v_{k} - v_{j} u_{k}}{\sqrt{\mbox{det} \abs{\mathbf{\Sigma}_{\mathbf{X}_{N} \mathbf{X}_{N}}}}}.
\end{align}
By virtue of \eqref{equation-39} and \eqref{equation-40}
\begin{align} \label{equation-41}
(\mathscr{E} (a_{j} b_{k} - b_{j} a_{k} \mid \mathbf{X}_{N} = \mathbf{K}))_{\substack{1 \leq j \leq N \\ 1 \leq k \leq N}}
 &= \frac{(K_{1}^{2} + K_{2}^{2}) (v_{j} u_{k} - u_{j} v_{k})}{\mbox{det} \abs{\mathbf{\Sigma}_{\mathbf{X}_{N} \mathbf{X}_{N}}}} \nonumber \\ &\quad + \frac{2 (u_{j} v_{k} - v_{j} u_{k})}{\sqrt{\mbox{det} \abs{\mathbf{\Sigma}_{\mathbf{X}_{N} \mathbf{X}_{N}}}}}.
\end{align}

It remains to evaluate the conditional expectation of $\mbox{det} \nabla \mathbf{X}_{N}$. From \eqref{equation-2}, \eqref{equation-35}, \eqref{equation-38}, and \eqref{equation-41}, it emerges from an arduous calculation that
\begin{equation*}
\begin{split}
&\mathscr{E} (\mbox{det} \nabla \mathbf{X}_{N} \mid \mathbf{X}_{N} = \mathbf{K}) \\
 &\hspace{25pt} = \sum_{j = 0}^{N} \mathscr{E} (a_{j}^{2} + b_{j}^{2} \mid \mathbf{X}_{N} = \mathbf{K}) \left\{\left(\frac{\partial u_{j}}{\partial x}\right)^{2} + \left(\frac{\partial v_{j}}{\partial x}\right)^{2}\right\} \\ &\hspace{35pt} + \sum_{j = 0}^{N} \sum_{\substack{k = 0 \\ k \neq j}}^{N} \left\{\mathscr{E} (a_{j} a_{k} + b_{j} b_{k} \mid \mathbf{X}_{N} = \mathbf{K}) \left(\frac{\partial u_{j}}{\partial x} \frac{\partial u_{k}}{\partial x} + \frac{\partial v_{j}}{\partial x} \frac{\partial v_{k}}{\partial x}\right)\right\} \\ &\hspace{35pt} + \sum_{j = 0}^{N} \sum_{\substack{k = 0 \\ k \neq j}}^{N} \left\{\mathscr{E} (a_{j} b_{k} - b_{j} a_{k} \mid \mathbf{X}_{N} = \mathbf{K}) \left(\frac{\partial v_{j}}{\partial x} \frac{\partial u_{k}}{\partial x} - \frac{\partial u_{j}}{\partial x} \frac{\partial v_{k}}{\partial x}\right)\right\} \\
 &\hspace{25pt} = 2 \sum_{j = 0}^{N} \abs{f_{j}^{\prime} (z)}^{2} - \left(\frac{2}{\sqrt{\mbox{det} \abs{\mathbf{\Sigma}_{\mathbf{X}_{N} \mathbf{X}_{N}}}}} - \frac{K_{1}^{2} + K_{2}^{2}}{\mbox{det} \abs{\mathbf{\Sigma}_{\mathbf{X}_{N} \mathbf{X}_{N}}}}\right) \\ &\hspace{115pt} \times \sum_{j = 0}^{N} \sum_{k = 0}^{N} \left\{(u_{j} u_{k} + v_{j} v_{k}) \left(\frac{\partial u_{j}}{\partial x} \frac{\partial u_{k}}{\partial x} + \frac{\partial v_{j}}{\partial x} \frac{\partial v_{k}}{\partial x}\right)\right. \\ &\hspace{155pt} \left. + (v_{j} u_{k} - u_{j} v_{k}) \left(\frac{\partial v_{j}}{\partial x} \frac{\partial u_{k}}{\partial x} - \frac{\partial u_{j}}{\partial x} \frac{\partial v_{k}}{\partial x}\right)\right\}.
\end{split}
\end{equation*}
An uninspired calculation shows that the double sum on the extreme right side can be reduced to
\begin{equation*}
\begin{split}
&\sum_{j = 0}^{N} \sum_{k = 0}^{N} \left(u_{j} \frac{\partial u_{j}}{\partial x} + v_{j} \frac{\partial v_{j}}{\partial x}\right) \left(u_{k} \frac{\partial u_{k}}{\partial x} + v_{k} \frac{\partial v_{k}}{\partial x}\right) \\ &\hspace{55pt} + \sum_{j = 0}^{N} \sum_{k = 0}^{N} \left(u_{j} \frac{\partial v_{j}}{\partial x} - v_{j} \frac{\partial u_{j}}{\partial x}\right) \left(u_{k} \frac{\partial v_{k}}{\partial x} - v_{k} \frac{\partial u_{k}}{\partial x}\right) \\
 &\hspace{45pt} = \left\{\sum_{j = 0}^{N} \left(u_{j} \frac{\partial u_{j}}{\partial x} + v_{j} \frac{\partial v_{j}}{\partial x}\right)\right\}^{2} + \left\{\sum_{j = 0}^{N} \left(u_{j} \frac{\partial v_{j}}{\partial x} - v_{j} \frac{\partial u_{j}}{\partial x}\right)\right\}^{2} \\
 &\hspace{45pt} = \left|\sum_{j = 0}^{N} \left\{\left(u_{j} \frac{\partial u_{j}}{\partial x} + v_{j} \frac{\partial v_{j}}{\partial x}\right) + i \left(u_{j} \frac{\partial v_{j}}{\partial x} - v_{j} \frac{\partial u_{j}}{\partial x}\right)\right\}\right|^{2} \\
 &\hspace{45pt} = \left|\sum_{j = 0}^{N} (u_{j} - i v_{j}) \left(\frac{\partial u_{j}}{\partial x} + i \frac{\partial v_{j}}{\partial x}\right)\right|^{2}
 = \left|\sum_{j = 0}^{N} \overline{f_{j} (z)} f_{j}^{\prime} (z)\right|^{2}.
\end{split}
\end{equation*}
Thus,
\begin{equation} \label{equation-42}
\begin{split}
&\mathscr{E} (\mbox{det} \nabla \mathbf{X}_{N} \mid \mathbf{X}_{N} = \mathbf{K}) \\
 &\quad = \frac{2}{\sqrt{\mbox{det} \abs{\mathbf{\Sigma}_{\mathbf{X}_{N} \mathbf{X}_{N}}}}} \left(\sqrt{\mbox{det} \abs{\mathbf{\Sigma}_{\mathbf{X}_{N} \mathbf{X}_{N}}}} \sum_{k = 0}^{N} \abs{f_{k}^{\prime} (z)}^{2} - \left|\sum_{j = 0}^{N} \overline{f_{j} (z)} f_{j}^{\prime} (z)\right|^{2}\right) \\ &\qquad + \frac{K_{1}^{2} + K_{2}^{2}}{\mbox{det} \abs{\mathbf{\Sigma}_{\mathbf{X}_{N} \mathbf{X}_{N}}}} \left|\sum_{j = 0}^{N} \overline{f_{j} (z)} f_{j}^{\prime} (z)\right|^{2}.
\end{split}
\end{equation}
For definiteness, we recall from \cite[Chapter 10]{Graybill1983} (see \cite[Chapter 2]{Tong1990}, also) that the joint density of two random complex Gaussian variables $X_{1, N}$ and $X_{2, N}$ at the points $x_{1}$ and $x_{2}$, respectively, is equal to
\begin{equation*}
p_{x, y} (x_{1}, x_{2})
 = \frac{1}{2 \pi \sigma^{2} \sqrt{\mbox{det} \abs{\mathbf{\Sigma}_{\mathbf{X}_{N} \mathbf{X}_{N}}}}} \exp \left(-\frac{(x_{1} - \mathscr{E} x_{1})^{2} + (x_{2} - \mathscr{E} x_{2})^{2}}{2 \sigma^{2} \sqrt{\mbox{det} \abs{\mathbf{\Sigma}_{\mathbf{X}_{N} \mathbf{X}_{N}}}}}\right).
\end{equation*}
In our case, the conditions $\mathscr{E} x_{1} = \mathscr{E} x_{2} = 0$ and $\sigma^{2} = 1$ apply. Thus, we have
\begin{equation} \label{equation-43}
p_{x, y} (\mathbf{K}^{\prime})
 = \frac{1}{2 \pi \sqrt{\mbox{det} \abs{\mathbf{\Sigma}_{\mathbf{X}_{N} \mathbf{X}_{N}}}}} \exp \left(-\frac{K_{1}^{2} + K_{2}^{2}}{\sqrt{\mbox{det} \abs{\mathbf{\Sigma}_{\mathbf{X}_{N} \mathbf{X}_{N}}}}}\right).
\end{equation}
By virtue of \eqref{equation-3}, \eqref{equation-42}, and \eqref{equation-43}, upon simplifying and applying the formulas for the $B_{r, N} (z)$ in Theorem \ref{theorem-1}, the required result follows.

\medskip

\section{The asymptotic analysis}

It is well known and, for example, Faramand \cite{Farahmand1986} has shown that, for large values of $N$, the real zeros of random polynomials with real coefficients are clustered about $\pm 1$. (See, also, Bharucha-Reid and Sambandham's book  \cite{Bharucha-ReidSambandham1986}.) In order to understand better the behaviour of $h_{N, \mathbf{K}} (z)$ in Theorem \ref{theorem-1} as $N$ tends to infinity, we define special values of the $f_{j} (z)$. Indeed, we are restricted to the cases that the evaluation of sums in Theorem \ref{theorem-1} becomes analytically feasible. To exhibit the numerical behaviour of $h_{N, \mathbf{K}} (z)$ and the zeros of the random sum $S_{N} (z)$ for various values of $N$ numerically, we used the general computing environment Wolfram Mathematica\textsuperscript{\textregistered} version number 12.0.0.0 developed by Wolfram Research for the platform Mac OS X x86 (64-bit), which ran on the Apple Mac Pro (late 2013) with the 2.7 GHz 12-core Intel\textsuperscript{\textregistered} Xeon\textsuperscript{\textregistered} Processor E5-2697 v2.

The simplest example of random sums is when
\begin{equation} \label{equation-44}
f_{j} (z)
 = z^{j}.
\end{equation}
These resulting limits are best expressed in terms of the function
\begin{equation*}
B (z)
 = \frac{1}{1 - \abs{z}^{2}}.
\end{equation*}
Using the notation given in Theorem \ref{theorem-1}, we have
\begin{equation*}
B_{0, N} (z)
 = \sum_{j = 0}^{N} \abs{z}^{2 j}
 = (1 - \abs{z}^{2 N + 2}) B (z).
\end{equation*}
By repeated differentiation, we obtain
\begin{equation*}
z B_{1, N} (z)
 = \sum_{j = 0}^{N} j \abs{z}^{2 j}
 = (N \abs{z}^{2 N + 2} - (N + 1) \abs{z}^{2 N} + 1) \abs{z}^{2} B (z)^{2}
\end{equation*}
and
\begin{equation*}
\begin{split}
B_{2, N} (z)
 = \sum_{j = 0}^{N} j^{2} \abs{z}^{2 j - 2}
 &= (1 + \abs{z}^{2} - \abs{z}^{2 N} (N^{2} \abs{z}^{4} - (2 N^{2} + 2 N - 1) \abs{z}^{2} \\ &\quad + (N + 1)^{2})) B (z)^{3}.
\end{split}
\end{equation*}
Clearly, if $\abs{z} < 1$, then
\begin{equation*}
\lim_{N \to \infty} B_{0, N} (z)
 = B (z)
 \end{equation*}
 and
 \begin{equation*}
 \lim_{N \to \infty} z B_{1, N} (z)
   = \abs{z}^{2} B (z)^{2},
\end{equation*}
as well as
\begin{equation*}
\lim_{N \to \infty} B_{2, N} (z)
 = (1 + \abs{z}^{2}) B (z)^{3}.
\end{equation*}
The following follows from Theorem \ref{theorem-1}.

\begin{theorem} \label{theorem-4}
Let the $f_{j} (z)$ in the definition of the random sum $S_{N} (z)$ be given by \eqref{equation-44}. If $\abs{z} < 1$, then we have
\begin{equation*}
\lim_{N \to \infty} h_{N, \mathbf{K}} (z)
 = \frac{1}{\pi} e^{-(K_{1}^{2} + K_{2}^{2}) / 2 B (z)} B (z) \left\{B (z) + \left(\frac{K_{1}^{2} + K_{2}^{2}}{2}\right) \abs{z}^{2}\right\}.
\end{equation*}
Then for any vector $\mathbf{K}$ restricted to a circle of radius $K > 0$ we have
\begin{equation*}
\lim_{N \to \infty} h_{N, \mathbf{K}} (z)
 = \frac{1}{\pi} e^{-K^{2} / B (z)} B (z) \left(B (z) + K^{2} \abs{z}^{2}\right),
\end{equation*}
and for any vector $\mathbf{K}$ of equal elements $K = 0$ we have
\begin{equation*}
\lim_{N \to \infty} h_{N, \mathbf{0}} (z)
 = \frac{1}{\pi} B (z)^{2}.
\end{equation*}
\end{theorem}

We note that, in all the cases considered, the limiting value of $h_{N, \mathbf{K}} (z)$ has $\abs{z}^{4} - 2 \abs{z}^{2} + 1$ in its denominators. An exponential factor is present when $\mathbf{K}$ is nonzero.

If now $\abs{z} > 1$, then for all sufficiently large $N$ we can write
\begin{equation*}
B_{0, N} (z)
 \sim -\abs{z}^{2 N} B (z)
\end{equation*}
and
\begin{equation*}
z B_{1, N} (z)
 \sim \abs{z}^{2} B (z)^{2},
\end{equation*}
as well as
\begin{equation*}
B_{2, N} (z)
 \sim (1 + \abs{z}^{2} - N^{2} \abs{z}^{2 N} (\abs{z}^{4} - \abs{z}^{2} + 1)) B (z)^{3}.
\end{equation*}
The following follows from Theorem \ref{theorem-1}.

\begin{theorem} \label{theorem-5}
Let the $f_{j} (z)$ in the definition of the random sum $S_{N} (z)$ be given by \eqref{equation-44}. If $\abs{z} < 1$, then for all sufficiently large $N$ we have
\begin{equation*}
\begin{split}
h_{N, \mathbf{K}} (z)
 &\sim \frac{1}{\pi} e^{(K_{1}^{2} + K_{2}^{2}) / 2 \abs{z}^{2 N} B (z)} \\ &\quad \times B (z)^{2} \left\{N^{2} (\abs{z}^{4} - \abs{z}^{2} + 1) - \frac{\abs{z}^{2} + 1}{\abs{z}^{2 N}} - \frac{\abs{z}^{2}}{\abs{z}^{4 N}} \left(1 + \frac{K_{1}^{2} + K_{2}^{2}}{2 \abs{z}^{2 N} B (z)}\right)\right\}.
\end{split}
\end{equation*}
Then for any vector $\mathbf{K}$ restricted to a circle of radius $K > 0$
\begin{equation*}
\begin{split}
h_{N, \mathbf{K}} (z)
 &\sim \frac{1}{\pi} e^{K^{2} / \abs{z}^{2 N} B (z)} \\ &\quad \times B (z)^{2} \left\{N^{2} (\abs{z}^{4} - \abs{z}^{2} + 1) - \frac{\abs{z}^{2} + 1}{\abs{z}^{2 N}} - \frac{\abs{z}^{2}}{\abs{z}^{4 N}} \left(1 + \frac{K^{2}}{\abs{z}^{2 N} B (z)}\right)\right\},
\end{split}
\end{equation*}
and for any vector $\mathbf{K}$ of equal elements $K = 0$
\begin{equation*}
h_{N, \mathbf{K}} (z)
 \sim \frac{1}{\pi} B (z)^{2} \left(N^{2} (\abs{z}^{4} - \abs{z}^{2} + 1) - \frac{\abs{z}^{2} + 1}{\abs{z}^{2 N}} - \frac{\abs{z}^{2}}{\abs{z}^{4 N}}\right).
\end{equation*}
\end{theorem}

Using the appropriate power sum formulas in Theorem \ref{theorem-1}, we obtain the following.

\begin{theorem} \label{theorem-6}
Let the $f_{j} (z)$ in the definition of the random sum $S_{N} (z)$ be given by \eqref{equation-44}. If $z = \pm 1$, then we have
\begin{equation*}
h_{N, \mathbf{K}} (\pm 1)
 = \frac{1}{12 \pi} e^{-(K_{1}^{2} + K_{2}^{2}) / (2 N + 2)} \left\{2 N + N^{2} \left(1 + \frac{3 (K_{1}^{2} + K_{2}^{2})}{2 N + 2}\right)\right\}.
\end{equation*}
Then for any vector $\mathbf{K}$ restricted to a circle of radius $K > 0$
\begin{equation*}
h_{N, \mathbf{K}} (\pm 1)
 = \frac{1}{12 \pi} e^{-K^{2} / (N + 1)} \left\{2 N + N^{2} \left( 1 + \frac{3 K^{2}}{N + 1}\right)\right\},
\end{equation*}
and for any vector $\mathbf{K}$ of equal elements $K = 0$
\begin{equation*}
h_{N, \mathbf{0}} (\pm 1)
 = \frac{1}{12\pi} N (N + 2).
\end{equation*}
\end{theorem}

\begin{figure}[!ht]
\begin{center}$
\begin{array}{cc}
\includegraphics[width=57mm]{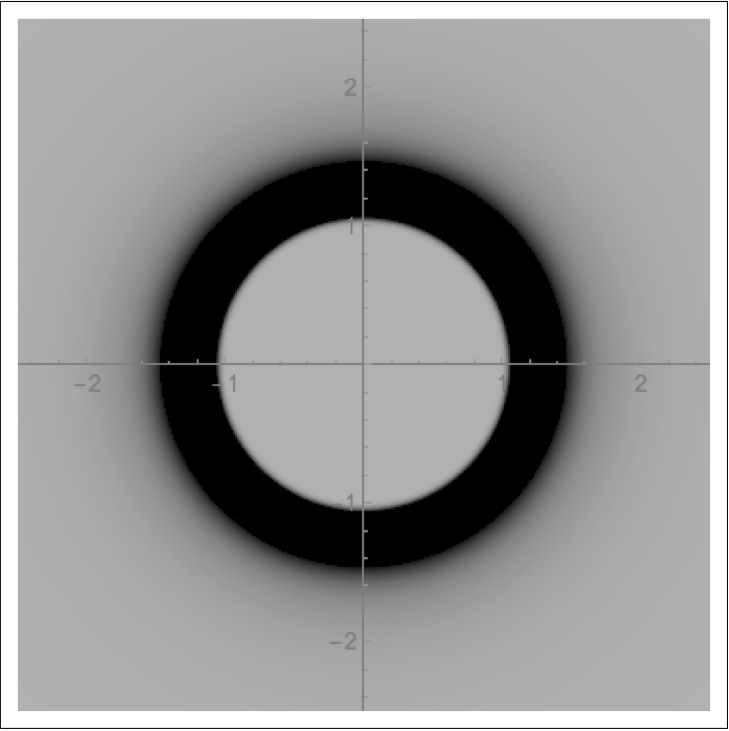} \hspace{25pt}
\includegraphics[width=57mm]{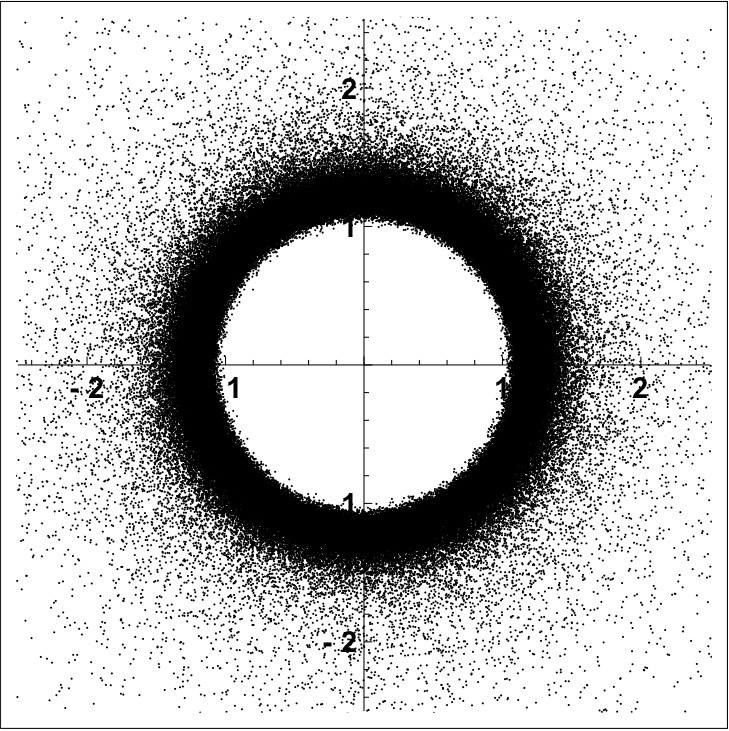}
\end{array}$
\end{center}
\caption{Twenty thousand random degree 10 polynomials for the random equation $S_{10} (z) = \eta_{0} + \eta_{1} z + \eta_{2} z^{2} + \ldots + \eta_{10} z^{10} = 10 + i 10$.}
\label{figure 1}
\end{figure}

In Figure \ref{figure 1}, the left-hand plot is a grey-scale image of ${h}_{N, \mathbf{K}}$ with $N = 10$ and $\mathbf{K} = (10, 10)^{\prime}$. The right-hand plot shows the zeros obtained by generating 20,000 random polynomials and explicitly finding the zeros of the random equation $S_{10} (z) = \mathbf{K}$. The zeros cluster near the unit circle, and the density function does not have mass concentrated on the real axis. There is no jump present near the real axis. For larger values of $\mathbf{K}$ the effect is more pronounced.

Also of interest are random Weyl polynomials in which
\begin{equation} \label{equation-45}
f_{j} (z)
 = \frac{z^{j}}{\sqrt{j!}},
\end{equation}
also studied by Farahmand and Jahangiri \cite{FarahmandJahangiri1998}, Littlewood and Offord \cite{LittlewoodOfford1938, LittlewoodOfford1939}, Offord \cite{Offord1995}, and Vanderbei \cite{Vanderbei2015}. For this case, the limiting forms of the various functions defining $h_{N, \mathbf{K}} (z)$ are computed by repeated differentiation. We have
\begin{equation*}
\lim_{N \to \infty} B_{0, N} (z)
 = \lim_{N \to \infty} \sum_{j = 0}^{N} \frac{\abs{z}^{2 j}}{j!}
 = e^{\abs{z}^{2}}
\end{equation*}
and
\begin{equation*}
\lim_{N \to \infty} B_{1, N} (z)
 = \lim_{N \to \infty} \sum_{j = 0}^{N} \frac{\abs{z}^{2 j}}{(j - 1)! z}
 = \overline{z} e^{\abs{z}^{2}},
\end{equation*}
as well as
\begin{equation*}
\lim_{N \to \infty} B_{2, N} (z)
 = \lim_{N \to \infty} \sum_{j = 0}^{N} \frac{j \abs{z}^{2 j - 2}}{(j - 1)!}
 = (\abs{z}^{2} + 1) e^{\abs{z}^{2}}.
\end{equation*}
Substituting these values into Theorem \ref{theorem-1}, we obtain the following. When $K = 0$, we recover Vanderbei's theorem \cite[Theorem 3.1]{Vanderbei2015}.

\begin{theorem} \label{theorem-7}
Let the $f_{j} (z)$ in the definition of the random sum $S_{N} (z)$ be given by \eqref{equation-45}. Then we have
\begin{equation*}
\lim_{N \to \infty} h_{N, \mathbf{K}} (z)
 = \frac{1}{\pi} e^{-(K_{1}^{2} + K_{2}^{2}) / 2 e^{\abs{z}^{2}}} \left(1 + \frac{(K_{1}^{2} + K_{2}^{2}) \abs{z}^{2}}{2 e^{\abs{z}^{2}}}\right).
\end{equation*}
Then for any vector $\mathbf{K}$ restricted to a circle of radius $K > 0$
\begin{equation*}
\lim_{N \to \infty} h_{N, \mathbf{K}} (z)
 = \frac{1}{\pi} e^{-K^{2} / e^{\abs{z}^{2}}} \left(1 + \frac{K^{2} \abs{z}^{2}}{e^{\abs{z}^{2}}}\right),
\end{equation*}
and for any vector $\mathbf{K}$ of equal elements $K = 0$
\begin{equation*}
\lim_{N \to \infty} h_{N, \mathbf{0}} (z)
 = \frac{1}{\pi}.
\end{equation*}
\end{theorem}

\begin{figure}[!ht]
\begin{center}$
\begin{array}{cc}
\includegraphics[width=57mm]{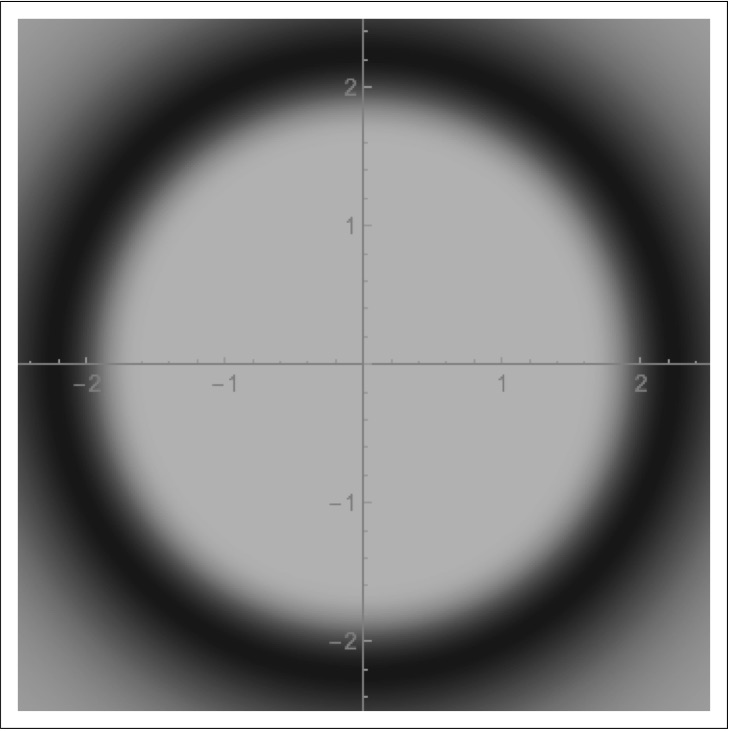} \hspace{25pt}
\includegraphics[width=57mm]{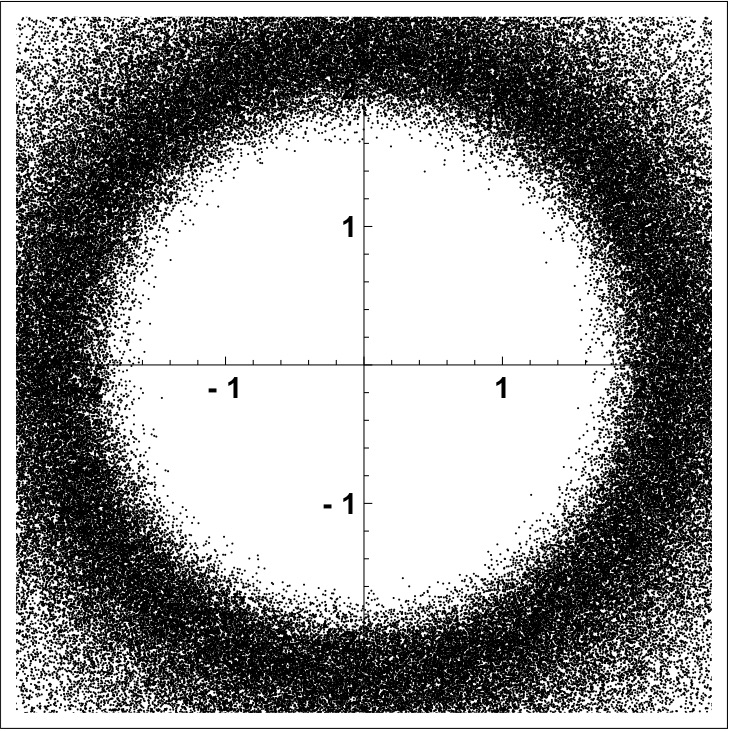}
\end{array}$
\end{center}
\caption{Twenty thousand random degree 10 Weyl polynomials for the random equation $S_{10} (z) = \eta_{0} + \eta_{1} z + \eta_{2} z^{2} / \sqrt{2!} + \ldots + \eta_{10} z^{10} / \sqrt{10!} = 10 + i 10$.}
\label{figure 2}
\end{figure}

We note that the distribution of the real zeros becomes uniform over the real number line. The complex zeros are much more uniformly distributed than was the case when the factor $1 / \sqrt{j!}$ was not present. The pictures in Figure \ref{figure 2} show the density function ${h}_{10, \mathbf{K}} (z)$ and the empirical distribution for 20,000 random sums when $f_{j} (z) = z^{j} / \sqrt{j!}$ and $\mathbf{K} = (10, 10)^{\prime}$, that is, random degree 10 Weyl polynomials. The behaviour of the density function and the empirical distribution for the random sums becomes very noticeable and intensified when $\mathbf{K}$ is increased.

\begin{figure}[!ht]
\begin{center}$
\begin{array}{cc}
\includegraphics[width=57mm]{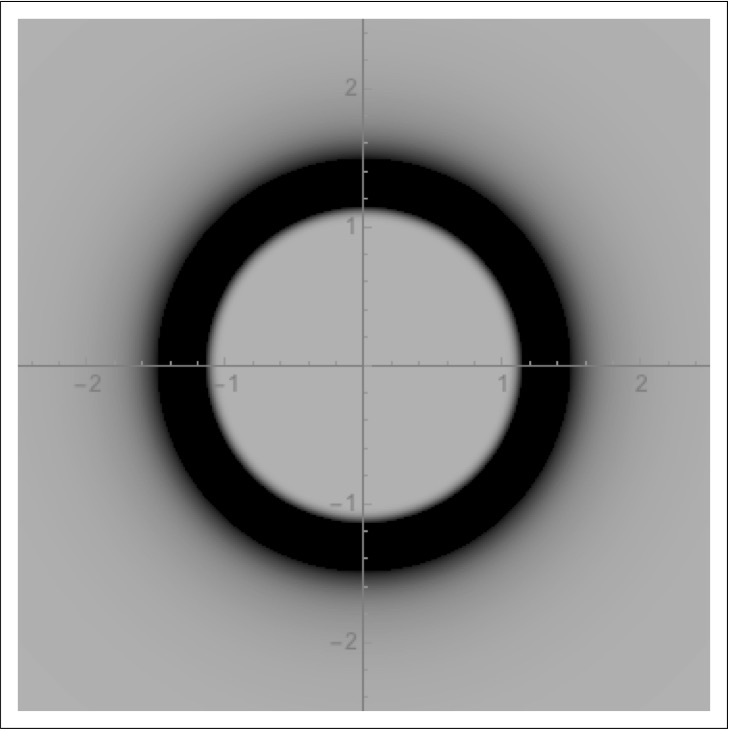} \hspace{25pt}
\includegraphics[width=57mm]{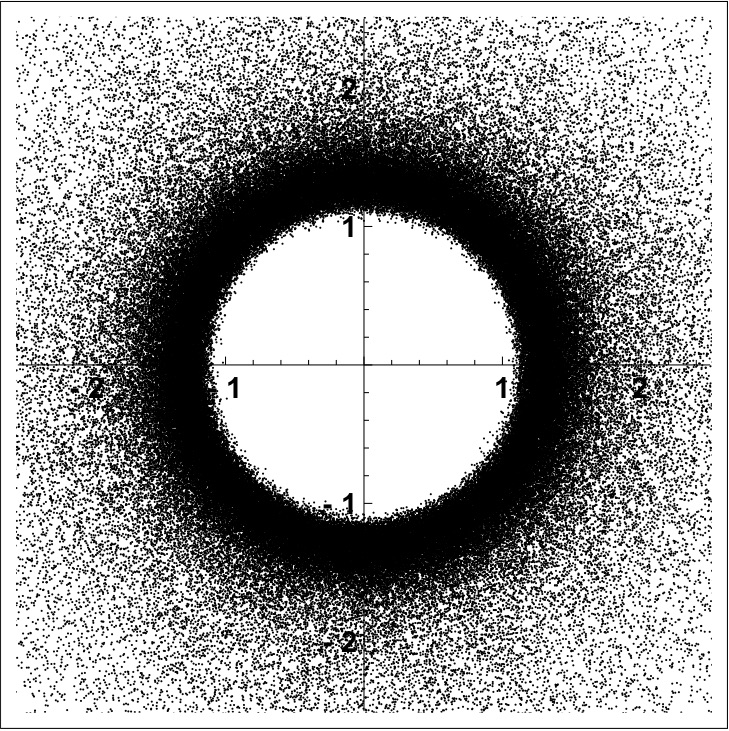}
\end{array}$
\end{center}
\caption{Twenty thousand random degree 10 root-binomial polynomials for the random equation $S_{10} (z) = \eta_{0} + \sqrt{5} \eta_{1} z + \sqrt{15} \eta_{2} z^{2} + \ldots + \sqrt{1 / 11} \eta_{10} z^{10} = 50 + i 50$.}
\label{figure 3}
\end{figure}

Next, let us assume that
\begin{equation} \label{equation-46}
f_{j} (z)
 = \sqrt{\binom{N}{j} \frac{1}{j + 1}} z^{j}.
\end{equation}
The random root-binomial polynomials were also studied by Farahmand and Jahangiri \cite{FarahmandJahangiri1998}, Littlewood and Offord \cite{LittlewoodOfford1938, LittlewoodOfford1939}, Offord \cite{Offord1995}, and Vanderbei \cite{Vanderbei2015}. The pictures in Figure \ref{figure 3} show the density function ${h}_{10, \mathbf{K}} (z)$ and the empirical distribution for 20,000 random degree 10 root-binomial polynomials with $\mathbf{K} = (50, 50)^{\prime}$.

By repeated differentiation, we can easily verify that
\begin{equation*}
B_{0, N} (z)
 = \sum_{j = 0}^{N} \binom{N}{j} \frac{\abs{z}^{2 j}}{j + 1}
 = \frac{(\abs{z}^{2} + 1)^{N + 1} - 1}{(N + 1) \abs{z}^{2}}
\end{equation*}
and
\begin{equation*}
B_{1, N} (z)
 = \sum_{j = 0}^{N} \binom{N}{j} \frac{j \abs{z}^{2 j}}{(j + 1) z}
 = \frac{(N \abs{z}^{2} - 1) (\abs{z}^{2} + 1)^{N} - 1}{(N + 1) z \abs{z}^{2}},
\end{equation*}
as well as
\begin{equation*}
B_{2, N} (z)
 = \sum_{j = 0}^{N} \binom{N}{j} \frac{j^{2} \abs{z}^{2 j - 2}}{j + 1}
 = \frac{(\abs{z}^{2} + 1)^{N} (\abs{z}^{2} (N^{2} \abs{z}^{2} - N + 1) + 1) - \abs{z}^{2} - 1}{(N + 1) \abs{z}^{4} (\abs{z}^{2} + 1)}.
\end{equation*}
If we now assume that $\abs{z} > 0$, that is, except at the origin, for all sufficiently large $N$ we can write
\begin{equation*}
B_{0, N} (z)
 \sim \frac{(\abs{z}^{2} + 1)^{N}}{N \abs{z}^{2}}
\end{equation*}
and
\begin{equation*}
B_{1, N} (z)
 \sim \frac{(\abs{z}^{2} + 1)^{N}}{z},
\end{equation*}
as well as
\begin{equation*}
B_{2, N} (z)
 \sim N (\abs{z}^{2} + 1)^{N - 1}.
\end{equation*}
We immediately get the following from Theorem \ref{theorem-1}.

\begin{theorem} \label{theorem-8}
Let the $f_{j} (z)$ in the definition of the random sum $S_{N} (z)$ be given by \eqref{equation-46}. Then for $\abs{z} > 0$ and all sufficiently large $N$ we have
\begin{equation*}
\begin{split}
h_{N, \mathbf{K}} (z)
 &\sim e^{-(K_{1}^{2} + K_{2}^{2}) N \abs{z}^{2} / 2 (\abs{z}^{2} + 1)^{N}} \frac{N^{2} \abs{z}^{4} ((K_{1}^{2} + K_{2}^{2}) N - 2 (\abs{z}^{2} + 1)^{N - 1})}{2 \pi (\abs{z}^{2} + 1)^{N}}.
\end{split}
\end{equation*}
Then for any vector $\mathbf{K}$ restricted to a circle of radius $K > 0$
\begin{equation*}
h_{N, \mathbf{K}} (z)
 \sim e^{-K^{2} N \abs{z}^{2} / (\abs{z}^{2} + 1)^{N}} \frac{N^{2} \abs{z}^{4} (K^{2} N - (\abs{z}^{2} + 1)^{N - 1})}{\pi (\abs{z}^{2} + 1)^{N}},
\end{equation*}
and for any vector $\mathbf{K}$ of equal elements $K = 0$
\begin{equation*}
h_{N, \mathbf{0}} (z)
 \sim -\frac{N^{2} \abs{z}^{4}}{\pi (\abs{z}^{2} + 1)}.
\end{equation*}
\end{theorem}

\begin{figure}[!ht]
\begin{center}$
\begin{array}{cc}
\includegraphics[width=57mm]{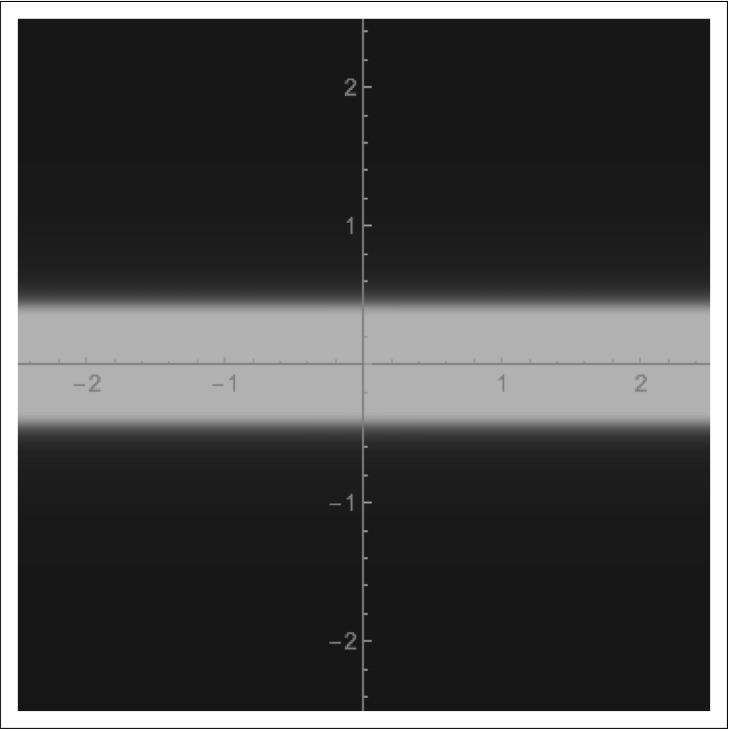} \hspace{25pt}
\includegraphics[width=57mm]{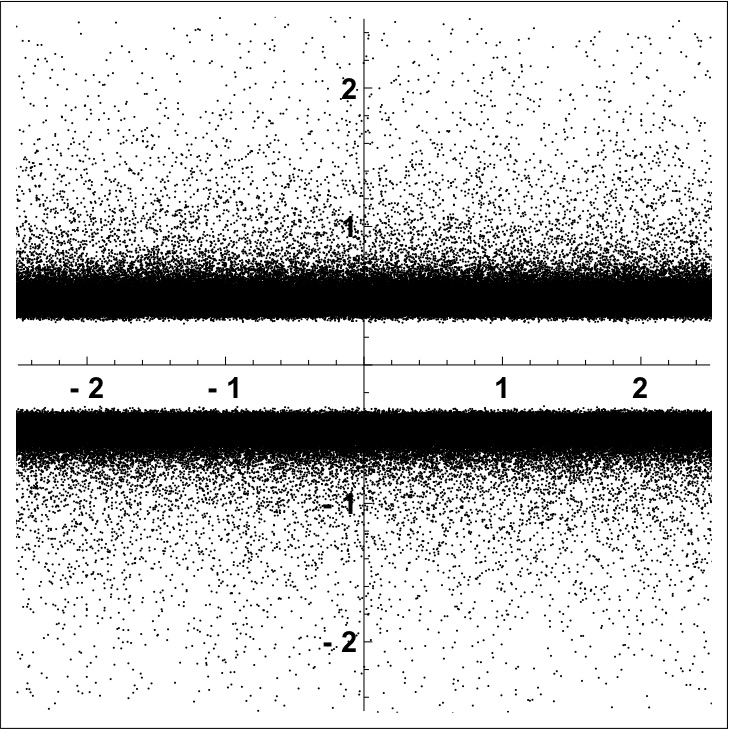}
\end{array}$
\end{center}
\caption{Twenty thousand random sums of the first 10 terms in a Fourier cosine series for the random equation $S_{10} = \eta_{0} + \eta_{1} \cos z + \eta_{2} \cos 2 z + \ldots + \eta_{10} \cos 10 z = 50 + i 50$.}
\label{figure 4}
\end{figure}

\begin{figure}[!ht]
\begin{center}$
\begin{array}{cc}
\includegraphics[width=57mm]{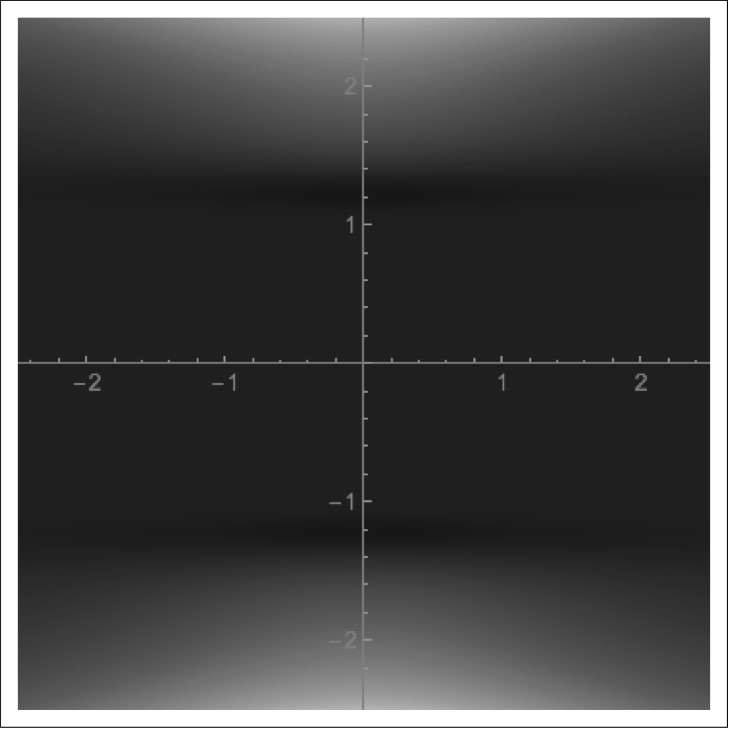} \hspace{25pt}
\includegraphics[width=57mm]{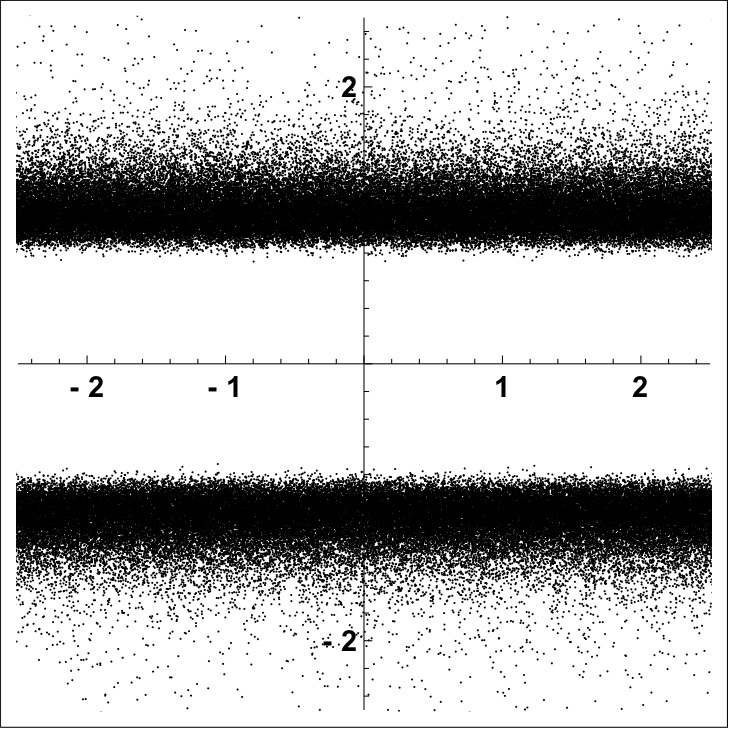}
\end{array}$
\end{center}
\caption{Twenty thousand random sums of the first 8 terms in a Fourier sine cosine series for the random equation $S_{8} (z) = \eta_{0} + \eta_{1} \sin z + \eta_{2} \cos z + \ldots + \eta_{7} \sin 4 z + \eta_{8} \cos 4 z = 50 + i 50$.}
\label{figure 5}
\end{figure}

Finally, we consider examples of random trigonometric sums. The behaviour of ${h}_{10, \mathbf{K}} (z)$ and the empirical distribution for a family of 20,000 random sums with $f_{j} (z) = \cos j z$ for $0 \leq j \leq 10$ and $\mathbf{K} = (50, 50)^{\prime}$ can be seen in Figure \ref{figure 4}. Since $\cos i y = \cosh y$, these random truncated Fourier cosine series are real-valued on both the real axis and the orthogonal imaginary axis. Figure \ref{figure 5} shows the corresponding behaviour for the random truncated Fourier sine and cosine series defined by
\begin{equation*}
f_{j} (z)
 = \left\{ \begin{array}{ll}
   \displaystyle \cos \left(\frac{j z}{2}\right) & \mbox{if $j$ is even,} \\ \addlinespace \addlinespace
   \displaystyle \sin \left(\frac{(j + 1) z}{2}\right) & \mbox{if $j$ is odd.}
\end{array}
\right.
\end{equation*}
The density functions for these examples for the case when $\mathbf{K} = (0, 0)^{\prime}$ have been obtained by Vanderbei \cite{Vanderbei2015}. We note that the density function for the random truncated Fourier sine and cosine series depends on the imaginary part of $z$ only.

\medskip

\section{The level crossings of random orthogonal polynomials}

We shall now consider the case when the $f_{j} (z)$ are either polynomials $p_{j} (z)$ orthogonal on the real line or polynomials $\varphi_{j} (z)$ orthogonal on the unit circle. These random orthogonal polynomials have real coefficients. These objects have been studied by many authors, including Das \cite{Das1971}, Lubinsky, Pritsker, and Xie \cite{LubinskyPritskerXie2016, LubinskyPritskerXie2018}, Yattselev and Yeager \cite{YattselevYeager2019}, and Yeager \cite{Yeager2016, Yeager2018}. In connection to these works, we are led to study the density functions for their complex zeros with respect to $\mathbf{K}$. First, we let $\alpha (z)$ denote a nondecreasing function with an infinite number of points of increase in the interval $[a, b]$. Assuming that moments of all orders exist, we say that a sequence of polynomials $\{p_{j} (z)\}_{j = 0}^{\infty}$, where the $p_{j} (z)$ have degree $n$, is orthogonal with respect to the distribution $\, d \alpha (z)$ if
\begin{equation*}
\int_{a}^{b} p_{j} (z) p_{k} (z) \, d \alpha (z)
 = h_{j} \delta_{k j}.
\end{equation*}
From \cite[Theorem 5.2.4]{AndrewsAskeyRoy1999}) and \cite[Theorem 3.2.2]{Szego1975}, the Christoffel--Darboux formula for orthogonal polynomials $p_{j} (z)$ on the real line can be stated as follows: Let $k_{j}$ be the highest coefficient of $p_{j} (z)$. Suppose that the $p_{j} (z)$ are normalized, so that
\begin{equation*}
h_{j}
 = \int_{a}^{b} p_{j} (z)^{2} \, d \alpha (z)
  = 1.
\end{equation*}
Then for complex variables $z$ and $w$
\begin{equation} \label{equation-47}
\sum_{j = 0}^{N} p_{j} (z) p_{j} (w)
 = \frac{k_{N}}{k_{N + 1}} \left(\frac{p_{N + 1} (z) p_{N} (w) - p_{N} (z) p_{N + 1} (w)}{z - w}\right).
\end{equation}
If $h_{j} \neq 1$, then \eqref{equation-47} takes the form
\begin{equation*}
\sum_{j = 0}^{N} \frac{p_{j} (z) p_{j} (w)}{h_{j}}
 = \frac{k_{N}}{k_{N + 1}} \left(\frac{p_{N + 1} (z) p_{N} (w) - p_{N} (z) p_{N + 1} (w)}{(z - w) h_{N}}\right).
\end{equation*}
(See \cite[Theorem 5.2.5]{AndrewsAskeyRoy1999}.) We have, also, the confluent form of \eqref{equation-47}, that is, when $z = w$. When $h_{j} = 1$, then
\begin{equation} \label{equation-48}
\sum_{j = 0}^{N} p_{j} (z)^{2}
 = \frac{k_{N}}{k_{N + 1}} (p_{N + 1}^{\prime} (z) p_{N} (z) - p_{N}^{\prime} (z) p_{N + 1} (z)).
\end{equation}

We find the representations for the $B_{r, N} (z)$. First, setting $w = \overline{z}$, so that $p_{j} (w) = p_{j} (\overline{z}) = \overline{p_{j} (z)}$, and $z - \overline{z} = 2 i \im (z)$ in \eqref{equation-47}, we obtain
\begin{equation} \label{equation-49}
B_{0, N} (z)
 = \sum_{j = 0}^{N} p_{j} (z) \overline{p_{j} (z)}
 = \frac{k_{N}}{k_{N + 1}} \left(\frac{\im (p_{N + 1} (z) \overline{p_{N} (z)})}{\im (z)}\right)
\end{equation}
Second, for $B_{1, N} (z)$ we first take the derivative of \eqref{equation-47} with respect to $w$ and then apply \eqref{equation-47} to get
\begin{align} \label{equation-50}
\sum_{j = 0}^{N} p_{j} (z) p_{j}^{\prime} (w)
 &= \frac{k_{N}}{k_{N + 1}} \left(\frac{p_{N + 1} (z) p_{N}^{\prime} (w) - p_{N} (z) p_{N + 1}^{\prime} (w)}{z - w} \right) \nonumber \\ &\quad + \frac{1}{z - w} \sum_{j = 0}^{N} p_{j} (z) p_{j} (w).
\end{align}
Then using $z$ in the place of $\overline{z}$ and putting $w = z$, so that $p_{j} (\overline{z}) = \overline{p_{j} (z)}$ and $p_{j}^{\prime} (w) = p_{j}^{\prime} (z)$, and $\overline{z} - z = -2 i \im (z)$ in \eqref{equation-50}, we obtain
\begin{align} \label{equation-51}
B_{1, N} (z)
 &= \sum_{j = 0}^{N} \overline{p_{j} (z)} p_{j}^{\prime} (z) \nonumber \\
 &= \frac{k_{N}}{k_{N + 1}} \left(\frac{\overline{p_{N} (z)} p_{N + 1}^{\prime} (z) - \overline{p_{N + 1} (z)} p_{N}^{\prime} (z)}{2 i \im (z)}\right) - \frac{B_{0, N} (z)}{2 i \im (z)}
\end{align}
Third, for $B_{2, N} (z)$ we first take the derivative of \eqref{equation-50} with respect to $z$ and then apply \eqref{equation-50} to receive
\begin{equation*}
\begin{split}
\sum_{j = 0}^{N} p_{j}^{\prime} (z) p_{j}^{\prime} (w)
 &= \frac{k_{N}}{k_{N + 1}} \left(\frac{p_{N + 1}^{\prime} (z) p_{N}^{\prime} (w) - p_{N}^{\prime} (z) p_{N + 1}^{\prime} (w)}{z - w}\right) \nonumber \\ &\quad - \frac{k_{N}}{k_{N + 1}} \left(\frac{p_{N + 1} (z) p_{N}^{\prime} (z) - p_{N} (z) p_{N + 1}^{\prime} (w)}{(z - w)^{2}}\right) \nonumber \\ &\quad + \frac{1}{z - w} \sum_{j = 0}^{N} p_{j}^{\prime} (z) p_{j} (w) - \frac{1}{(z - w)^{2}} \sum_{j = 0}^{N} p_{j} (z) p_{j} (w).
\end{split}
\end{equation*}
Again, from \eqref{equation-50}, since
\begin{equation*}
\begin{split}
&\frac{k_{N}}{k_{N + 1}} \left(\frac{p_{N + 1} (z) p_{N}^{\prime} (z) - p_{N} (z) p_{N + 1}^{\prime} (w)}{(z - w)^{2}}\right) \\
 &\hspace{55pt}= \frac{1}{z - w} \sum_{j = 0}^{N} p_{j} (z) p_{j}^{\prime} (w) - \frac{1}{(z - w)^{2}} \sum_{j = 0}^{N} p_{j} (z) p_{j} (w),
\end{split}
\end{equation*}
we have
\begin{equation*}
\begin{split}
\sum_{j = 0}^{N} p_{j}^{\prime} (z) p_{j}^{\prime} (w)
 &= \frac{k_{N}}{k_{N + 1}} \left(\frac{p_{N + 1}^{\prime} (z) p_{N}^{\prime} (w) - p_{N}^{\prime} (z) p_{N + 1}^{\prime} (w)}{z - w}\right) \\ &\quad - \frac{1}{z - w} \sum_{j = 0}^{N} p_{j} (z) p_{j}^{\prime} (w) + \frac{1}{z - w} \sum_{j = 0}^{N} p_{j}^{\prime} (z) p_{j} (w).
\end{split}
\end{equation*}
Thus,
\begin{align} \label{equation-52}
B_{2, N} (z)
 &= \sum_{j = 0}^{N} p_{j}^{\prime} (z) \overline{p_{j}^{\prime} (z)} \nonumber \\
 &= \frac{k_{N}}{k_{N + 1}} \left(\frac{\im (p_{N + 1}^{\prime} (z) \overline{p_{N}^{\prime} (z)})}{\im (z)}\right) - \frac{\overline{B_{1, N} (z)}}{2 i \im (z)} + \frac{B_{1, N} (z)}{2 i \im (z)}.
\end{align}

We apply \eqref{equation-49}, \eqref{equation-51}, and \eqref{equation-52} to Theorem \ref{theorem-1} to obtain the following formula for the density function $h_{N, \mathbf{K}} (z)$ for the complex zeros of random polynomials orthogonal on the real line.

\begin{theorem} \label{theorem-9}
Let the sequence $\{f_{j} (z)\}_{j = 0}^{N}$ in the definition of the random sum $S_{N} (z)$ be polynomials $p_{j} (z)$ orthogonal on the real line. For all integers $N > 0$ we have
\begin{equation*}
\begin{split}
&h_{N, \mathbf{K}} (z)
 = \frac{1}{\pi} \exp\left(-\frac{(K_{1}^{2} + K_{2}^{2}) k_{N + 1} \im (z)}{k_{N} \im (p_{N + 1} (z) \overline{p_{N} (z)})}\right) \\ &\hspace{15pt} \times \Vast\{\frac{\im (p_{N + 1}^{\prime} (z) \overline{p_{N}^{\prime} (z)})}{\im (p_{N + 1} (z) \overline{p_{N} (z)})} - \frac{\abs{\overline{p_{N} (z)} p_{N + 1}^{\prime} (z) - \overline{p_{N + 1} (z)} p_{N}^{\prime} (z)}^{2}}{4 \im (p_{N + 1} (z) \overline{p_{N} (z)})^{2}} + \frac{1}{4 \im (z)^{2}} \\ &\hspace{25pt} + \frac{(K_{1}^{2} + K_{2}^{2}) k_{N + 1}}{k_{N}} \vast(\frac{\abs{\overline{p_{N} (z)} p_{N + 1}^{\prime} (z) - \overline{p_{N + 1} (z)} p_{N}^{\prime} (z)}^{2} \im (z)}{8 \im (p_{N + 1} (z) \overline{p_{N} (z)})^{3}} \\ &\hspace{35pt} - \frac{\re (p_{N} (z) p_{N + 1}^{\prime} (z) - p_{N + 1} (z) \overline{p_{N}^{\prime} (z)})}{4 \im (p_{N + 1} (z) \overline{p_{N} (z)})^{2}} + \frac{1}{8 \im (z) \im (p_{N + 1} (z) \overline{p_{N} (z)})}\vast) \Vast\}.
\end{split}
\end{equation*}
When $\mathbf{K} = \mathbf{0}$, then
\begin{equation*}
\begin{split}
&h_{N, \mathbf{K}} (z) \\
& \hspace{15pt} = \frac{1}{\pi} \left(\frac{\im (p_{N + 1}^{\prime} (z) \overline{p_{N}^{\prime} (z)})}{\im (p_{N + 1} (z) \overline{p_{N} (z)})} - \frac{\abs{\overline{p_{N} (z)} p_{N + 1}^{\prime} (z) - \overline{p_{N + 1} (z)} p_{N}^{\prime} (z)}^{2}}{4 \im (p_{N + 1} (z) \overline{p_{N} (z)})^{2}} + \frac{1}{4 \im (z)^{2}}\right).
\end{split}
\end{equation*}
\end{theorem}

Next, the Christoffel--Darboux formula for a  sequence of polynomials $\{\varphi_{j} (z)\}_{j = 0}^{\infty}$ orthogonal on the unit circle is that, for complex variables $z$ and $w$ with $\overline{w} z \neq 1$,
\begin{equation} \label{equation-53}
\sum_{j = 0}^{N} \varphi_{j} (z) \overline{\varphi_{j} (w)}
 = \frac{\overline{\varphi_{N + 1}^{\ast} (w)} \varphi_{N + 1}^{\ast} (z) - \overline{\varphi_{N + 1} (w)} \varphi_{N + 1} (z)}{1 - \overline{w} z},
\end{equation}
where
\begin{equation} \label{equation-54}
\varphi_{N}^{\ast}
 = z^{N} \overline{\varphi_{N} \left(\frac{1}{\overline{z}}\right)}.
\end{equation}
As before, we find the representations for the $B_{r, N} (z)$. First, from \eqref{equation-53}
\begin{equation} \label{equation-55}
B_{0, N} (z)
 = \sum_{j = 0}^{N} \varphi_{j} (z) \overline{\varphi_{j} (z)}
 = \frac{\abs{\varphi_{N + 1}^{\ast} (z)}^{2} - \abs{\varphi_{N + 1} (z)}^{2}}{1 - \abs{z}^{2}}.
\end{equation}
Second, for $B_{1, N} (z)$ we first take the derivative of \eqref{equation-53} with respect to $\overline{w}$ and use $\overline{\varphi_{N + 1}^{\ast} (w)} = \varphi_{N + 1}^{\ast} (\overline{w})$ and $\overline{\varphi_{N + 1} (w)} = \varphi_{N + 1} (\overline{w})$ to obtain
\begin{equation} \label{equation-56}
\begin{split}
\sum_{j = 0}^{N} \varphi_{j} (z) \overline{\varphi_{j}^{\prime} (w)}
 &= \frac{\overline{\varphi_{N + 1}^{\ast \prime} (w)} \varphi_{N + 1}^{\ast} (z) - \overline{\varphi_{N + 1}^{\prime} (w)} \varphi_{N + 1} (z)}{1 - \overline{w} z} \\ &\quad + \frac{z (\overline{\varphi_{N + 1}^{\ast} (w)} \varphi_{N + 1}^{\ast} (z) - \overline{\varphi_{N + 1} (w)} \varphi_{N + 1} (z))}{(1 - \overline{w} z)^{2}}
\end{split}
\end{equation}
Then putting $w = z$ in \eqref{equation-56} and applying \eqref{equation-55} we obtain
\begin{align} \label{equation-57}
\overline{B_{1, N} (z)}
 &= \sum_{j = 0}^{N} \varphi_{j} (z) \overline{\varphi_{j}^{\prime} (z)} \nonumber \\
 &= \frac{\overline{\varphi_{N + 1}^{\ast \prime} (z)} \varphi_{N + 1}^{\ast} (z) - \overline{\varphi_{N + 1}^{\prime} (z)} \varphi_{N + 1} (z)}{1 - \abs{z}^{2}} + \frac{z B_{0, N} (z)}{1 - \abs{z}^{2}}.
\end{align}
Third, for $B_{2, N} (z)$ we first take the derivative of \eqref{equation-56} with respect to $z$ to obtain
\begin{equation} \label{equation-58}
\begin{split}
\sum_{j = 0}^{N} \varphi^{\prime} (z) \overline{\varphi_{j}^{\prime} (w)}
 &= \frac{\overline{\varphi_{N + 1}^{\ast \prime} (w)} \varphi_{N + 1}^{\ast \prime} (z) - \overline{\varphi_{N + 1}^{\prime} (w)} \varphi_{N + 1}^{\prime} (z)}{1 - \overline{w} z} \\ &\quad + \frac{\overline{w} (\overline{\varphi_{N + 1}^{\ast \prime} (w)} \varphi_{N + 1}^{\ast} (z) - \overline{\varphi_{N + 1} (w)} \varphi_{N + 1} (z))}{(1 - \overline{w} z)^{2}} \\ &\quad + \frac{z (\overline{\varphi_{N + 1}^{\ast} (w)} \varphi_{N + 1}^{\ast} (z) - \overline{\varphi_{N + 1} (w)} \varphi_{N + 1}^{\prime} (z))}{(1 - \overline{w} z)^{2}} \\ &\quad + \frac{\overline{\varphi_{N + 1}^{\ast} (w)} \varphi_{N + 1}^{\ast} (z) - \overline{\varphi_{N + 1} (w)} \varphi_{N + 1} (z)}{(1 - \overline{w} z)^{2}} \\ &\quad + \frac{2 \overline{w} z (\overline{\varphi_{N + 1}^{\ast} (w)} \varphi_{N + 1}^{\ast} (z) - \overline{\varphi_{N + 1} (w)} \varphi_{N + 1} (z))}{(1 - \overline{w} z)^{3}}.
\end{split}
\end{equation}
Thus, putting $w = z$ in \eqref{equation-58}, applying \eqref{equation-55} and \eqref{equation-57}, and noting that
\begin{equation*}
\frac{\overline{z} (\overline{\varphi_{N + 1}^{\ast \prime} (z)} \varphi_{N + 1}^{\ast} (z) - \overline{\varphi_{N + 1} (z)} \varphi_{N + 1} (z))}{(1 - \abs{z}^{2})^{2}}
 = \frac{\overline{z B_{1, N} (z)}}{1 - \abs{z}^{2}} - \frac{\overline{z} B_{0, N} (z)}{1 - \abs{z}^{2}}
\end{equation*}
and
\begin{equation*}
\frac{z (\overline{\varphi_{N + 1}^{\ast} (z)} \varphi_{N + 1}^{\ast} (z) - \overline{\varphi_{N + 1} (z)} \varphi_{N + 1}^{\prime} (z))}{(1 - \abs{z}^{2})^{2}}
 = \frac{z B_{0, N} (z)}{1 - \abs{z}^{2}} - \frac{\abs{z}^{2} B_{0, N} (z)}{(1 - \abs{z}^{2})^{2}}
\end{equation*}
we obtain
\begin{align} \label{equation-59}
B_{2, N} (z)
 &= \sum_{j = 0}^{N} \varphi_{j}^{\prime} (z) \overline{\varphi_{j}^{\prime} (z)} \nonumber \\
 &= \frac{\abs{\varphi_{N + 1}^{\ast \prime} (z)}^{2} - \abs{\varphi_{N + 1}^{\prime} (z)}^{2}}{1 - \abs{z}^{2}} + \frac{2 \re (z B_{1, N} (z))}{1 - \abs{z}^{2}} + \frac{B_{0, N} (z)}{1 - \abs{z}^{2}}.
\end{align}
To facilitate the derivation of the density function $h_{N, \mathbf{K}} (z)$ for the complex zeros of random polynomials orthogonal on the unit circle, we use the fact that the formula for $h_{N, \mathbf{K}} (z)$ in Theorem \ref{theorem-1} contains the quotients $(B_{0, N} (z) B_{2, N} (z) - \abs{B_{1, N} (z)}^{2}) / B_{0, N} (z)^{2}$ and $\abs{B_{1, N} (z)}^{2} / B_{0, N} (z)^{3}$. We treat these quotients in turn.

From \eqref{equation-55} and \eqref{equation-59}
\begin{equation} \label{equation-60}
\begin{split}
B_{0, N} (z) B_{2, N} (z)
 &= \frac{\abs{\varphi_{N + 1}^{\ast \prime} (z) \varphi_{N + 1}^{\ast} (z)}^{2}}{(1 - \abs{z}^{2})^{2}} + \frac{\abs{\varphi_{N + 1}^{\prime} (z) \varphi_{N + 1} (z)}^{2}}{(1 - \abs{z}^{2})^{2}} \\ &\quad - \frac{\abs{\varphi_{N + 1}^{\ast \prime} (z) \varphi_{N + 1} (z)}^{2}}{(1 - \abs{z}^{2})^{2}} - \frac{\abs{\varphi_{N + 1}^{\ast} (z) \varphi_{N + 1} (z)}^{2}}{(1 - \abs{z}^{2})^{2}} \\ &\quad + \frac{2 B_{0, N} \re (z B_{1, N} (z))}{1 - \abs{z}^{2}} + \frac{B_{0, N} (z)^{2}}{1 - \abs{z}^{2}}.
\end{split}
\end{equation}
From \eqref{equation-57}
\begin{align} \label{equation-61}
\abs{B_{1, N} (z)}^{2}
 &= \overline{B_{1, N} (z)} B_{1, N} (z) \nonumber \\
 &= \frac{\abs{\varphi_{N + 1}^{\ast \prime} \varphi_{N + 1}^{\ast} (z)}^{2}}{(1 - \abs{z}^{2})^{2}} + \frac{\abs{\varphi_{N + 1}^{\prime} (z) \varphi_{N + 1} (z)}^{2}}{(1 - \abs{z}^{2})^{2}} \nonumber \\ &\quad - \frac{2 \re (\overline{\varphi_{N + 1}^{\ast \prime} (z)} \varphi_{N + 1}^{\ast} (z) \varphi_{N + 1}^{\prime} (z) \overline{\varphi_{N + 1} (z)})}{(1 - \abs{z}^{2})^{2}} \nonumber \\ &\quad + \frac{2 B_{0, N} (z) \re (z B_{1, N} (z))}{1 - \abs{z}^{2}} - \frac{\abs{z}^{2} B_{0, N} (z)^{2}}{(1 - \abs{z}^{2})^{2}}.
\end{align}
Thus, from \eqref{equation-60} and \eqref{equation-61}
\begin{equation} \label{equation-62}
\begin{split}
B_{0, N} (z) B_{2, N} (z) - \abs{B_{1, N} (z)}^{2}
  &= \frac{(\abs{\varphi_{N + 1}^{\ast} (z)}^{2} - \abs{\varphi_{N + 1} (z)}^{2})^{2}}{(1 - \abs{z}^{2})^{4}} \\ &\quad - \frac{\abs{\varphi_{N + 1}^{\ast} (z) \varphi_{N + 1}^{\prime} (z) - \varphi_{N + 1}^{\ast \prime} (z) \varphi_{N + 1} (z)}^{2}}{(1 - \abs{z}^{2})^{2}}.
\end{split}
\end{equation}
From \eqref{equation-55} and \eqref{equation-61}
\begin{equation} \label{equation-63}
\begin{split}
\frac{\abs{B_{1, N} (z)}^{2}}{B_{0, N} (z)^{3}}
 &= \frac{(1 - \abs{z}^{2}) \abs{\overline{\varphi_{N + 1}^{\ast \prime} (z)} \varphi_{N + 1}^{\ast} (z) - \overline{\varphi_{N + 1}^{\prime} (z)} \varphi_{N + 1} (z)}^{2}}{(\abs{\varphi_{N + 1}^{\ast} (z)}^{2} - \abs{\varphi_{N + 1} (z)}^{2})^{3}} \\ &\quad + \frac{2 \re (z (\varphi_{N + 1}^{\ast \prime} (z) \overline{\varphi_{N + 1}^{\ast} (z)} - \varphi_{N + 1}^{\prime} (z) \overline{\varphi_{N + 1} (z)}))}{(\abs{\varphi_{N + 1}^{\ast} (z)}^{2} - \abs{\varphi_{N + 1} (z)}^{2})^{2}} \\ &\quad + \frac{\abs{z}^{2}}{(1 - \abs{z}^{2}) (\abs{\varphi_{N + 1}^{\ast} (z)}^{2} - \abs{\varphi_{N + 1} (z)}^{2})}.
\end{split}
\end{equation}
We deduce from Theorem \ref{theorem-1}, \eqref{equation-55}, \eqref{equation-62}, and \eqref{equation-63} the following result.

\begin{theorem} \label{theorem-10}
Let the sequence $\{f_{j} (z)\}_{j = 0}^{N}$ in the definition of the random sum $S_{N} (z)$ be polynomials $\varphi_{j} (z)$ orthogonal on the unit circle. Let $\varphi^{\ast}$ be given by \eqref{equation-54}. For all integers $N > 0$ we have
\begin{equation*}
\begin{split}
h_{N, \mathbf{K}} (z)
 &= \frac{1}{\pi} \exp \left\{-\left(\frac{K_{1}^{2} + K_{2}^{2}}{2}\right) \frac{1 - \abs{z}^{2}}{\abs{\varphi_{N + 1}^{\ast} (z)}^{2} - \abs{\varphi_{N + 1} (z)}^{2}}\right\} \\ &\hspace{15pt} \times \Vast\{\frac{1}{(1 - \abs{z}^{2})^{2}} - \frac{\abs{\varphi_{N + 1}^{\ast} (z) \varphi_{N + 1}^{\prime} (z) - \varphi_{N + 1}^{\ast \prime} (z) \varphi_{N + 1} (z)}^{2}}{(\abs{\varphi_{N + 1}^{\ast} (z)}^{2} - \abs{\varphi_{N + 1} (z)}^{2})^{2}} \\ &\hspace{25pt} + \left(\frac{K_{1}^{2} + K_{2}^{2}}{2}\right) \vast(\frac{(1 - \abs{z}^{2}) \abs{\overline{\varphi_{N + 1}^{\ast \prime} (z)} \varphi_{N + 1}^{\ast} (z) - \overline{\varphi_{N + 1}^{\prime} (z)} \varphi_{N + 1} (z)}^{2}}{(\abs{\varphi_{N + 1}^{\ast}}^{2} - \abs{\varphi_{N + 1} (z)}^{2})^{3}} \\ &\hspace{75pt} + \frac{2 \re (z (\varphi_{N + 1}^{\ast \prime} (z) \overline{\varphi_{N + 1}^{\ast} (z)} - \varphi_{N + 1}^{\prime} (z) \overline{\varphi_{N + 1} (z)}))}{(\abs{\varphi_{N + 1}^{\ast} (z)}^{2} - \abs{\varphi_{N + 1} (z)}^{2})^{2}} \\ &\hspace{125pt} + \frac{\abs{z}^{2}}{(1 - \abs{z}^{2}) (\abs{\varphi_{N + 1}^{\ast} (z)}^{2} - \abs{\varphi_{N + 1} (z)}^{2}))}\vast)\Vast\}.
\end{split}
\end{equation*}
When $\mathbf{K} = \mathbf{0}$, then
\begin{equation*}
h_{N, \mathbf{K}} (z) \\
 = \frac{1}{\pi} \left(\frac{1}{(1 - \abs{z}^{2})^{2}} - \frac{\abs{\varphi_{N + 1}^{\ast} (z) \varphi_{N + 1}^{\prime} (z) - \varphi_{N + 1}^{\ast \prime} (z) \varphi_{N + 1} (z)}^{2}}{(\abs{\varphi_{N + 1}^{\ast} (z)}^{2} - \abs{\varphi_{N + 1} (z)}^{2})^{2}}\right).
\end{equation*}
\end{theorem}

\medskip

\section{Final comments and future research}
The method introduced by Shepp and Vanderbei \cite{SheppVanderbei1995} based on Cauchy's argument principle could be applied in many circumstances. The method could be modified to produce the number of $\mathbf{K}$ complex level crossings. Furthermore, working with real coefficients, in fact, makes the analysis more complicated. This will be addressed in a future paper.

\medskip

\bibliographystyle{amsplain}

\end{document}